\documentclass[12pt]{amsart}

\input bookman.sty
\boldmath

\usepackage{helvet}

\usepackage{graphics}
\usepackage{amssymb}
\usepackage{amsxtra}
\usepackage{amsmath}
\usepackage{mathrsfs}

\usepackage[arrow, matrix]{xy}
\xyoption{frame}

\newfont{\bff}{cmbx10  scaled 1000}

\theoremstyle{plain}

\newtheorem{theo}{Theorem}[section]
\newtheorem{lemm}[theo]{Lemma}
\newtheorem{prop}[theo]{Proposition}
\newtheorem{coro}[theo]{Corollary}

\theoremstyle{definition}

\newtheorem{defi}[theo]{Definition}
\newtheorem{rema}[theo]{Remark}

\newfont{\rmm}{cmr10 scaled 1000}

\newfont{\itt}{cmsl10 scaled 1000}

\newfont{\rM}{cmr10 scaled 1700}

\newcounter{lemma}[section]

\newcounter{tempcounter}

\newcommand{\lb}{\label}

\newcommand{\rrf}[1]{(\ref{#1})}

\begin{document}


\renewcommand{\a}{\alpha}
\renewcommand{\b}{\beta}
\newcommand{\g}{\gamma}
\renewcommand{\d}{\delta}
\newcommand{\e}{\epsilon}
\newcommand{\ve}{\varepsilon}
\newcommand{\z}{\zeta}
\renewcommand{\t}{\theta}
\renewcommand{\l}{\lambda}
\renewcommand{\k}{\varkappa}
\newcommand{\m}{\mu}
\newcommand{\n}{\nu}
\renewcommand{\r}{\rho}
\newcommand{\vr}{\varrho}
\newcommand{\s}{\sigma}
\newcommand{\vp}{\varphi}
\renewcommand{\o}{\omega}

\renewcommand{\Re}{\text{\rm Re }}

\newcommand{\G}{\Gamma}
\newcommand{\D}{\Delta}
\newcommand{\T}{\Theta}
\renewcommand{\L}{\Lambda}
\renewcommand{\P}{\Pi}
\newcommand{\Si}{\Sigma}
\renewcommand{\O}{\Omega}

\newcommand{\Up}{\Upsilon}

\renewcommand{\AA}{{\mathcal A}}
\newcommand{\BB}{{\mathcal B}}
\newcommand{\CC}{{\mathcal C}}
\newcommand{\DD}{{\mathcal D}}
\newcommand{\EE}{{\mathcal E}}
\newcommand{\FF}{{\mathcal F}}
\newcommand{\GG}{{\mathcal G}}
\newcommand{\HH}{{\mathcal H}}
\newcommand{\II}{{\mathcal I}}
\newcommand{\JJ}{{\mathcal J}}
\newcommand{\KK}{{\mathcal K}}
\newcommand{\LL}{{\mathcal L}}
\newcommand{\MM}{{\mathcal M}}
\newcommand{\NN}{{\mathcal N}}
\newcommand{\OO}{{\mathcal O}}
\newcommand{\PP}{{\mathcal P}}
\newcommand{\QQ}{{\mathcal Q}}
\newcommand{\RR}{{\mathcal R}}
\renewcommand{\SS}{{\mathcal S}}
\newcommand{\TT}{{\mathcal T}}
\newcommand{\UU}{{\mathcal U}}
\newcommand{\VV}{{\mathcal V}}
\newcommand{\WW}{{\mathcal W}}
\newcommand{\XX}{{\mathcal X}}
\newcommand{\YY}{{\mathcal Y}}
\newcommand{\ZZ}{{\mathcal Z}}

\renewcommand{\aa}{{\mathbb{A}}}
\newcommand{\bb}{{\mathbb{B}}}
\newcommand{\cc}{{\mathbb{C}}}
\newcommand{\dd}{{\mathbb{D}}}
\newcommand{\ee}{{\mathbb{E}}}
\newcommand{\ff}{{\mathbb{F}}}
\renewcommand{\gg}{{\mathbb{G}}}
\newcommand{\hh}{{\mathbb{H}}}
\newcommand{\ii}{{\mathbb{I}}}
\newcommand{\jj}{{\mathbb{J}}}
\newcommand{\kk}{{\mathbb{K}}}
\renewcommand{\ll}{{\mathbb{L}}}
\newcommand{\mm}{{\mathbb{M}}}
\newcommand{\nn}{{\mathbb{N}}}
\newcommand{\oo}{{\mathbb{O}}}
\newcommand{\pp}{{\mathbb{P}}}
\newcommand{\qq}{{\mathbb{Q}}}
\newcommand{\rr}{{\mathbb{R}}}
\renewcommand{\ss}{{\mathbb{S}}}
\newcommand{\ttt}{{\mathbb{T}}}
\newcommand{\uu}{{\mathbb{U}}}
\newcommand{\vv}{{\mathbb{V}}}
\newcommand{\ww}{{\mathbb{W}}}
\newcommand{\xx}{{\mathbb{X}}}
\newcommand{\yy}{{\mathbb{Y}}}
\newcommand{\zz}{{\mathbb{Z}}}

\newcommand{\AAA}{{\mathbf{A}}}
\newcommand{\BBB}{{\mathbf{B}} }
\newcommand{\CCC}{{\mathbf{C}} }
\newcommand{\DDD}{{\mathbf{D}} }
\newcommand{\EEE}{{\mathbf{E}} }
\newcommand{\FFF}{{\mathbf{F}} }
\newcommand{\GGG}{{\mathbf{G}}}
\newcommand{\HHH}{{\mathbf{H}}}
\newcommand{\III}{{\mathbf{I}}}
\newcommand{\JJJ}{{\mathbf{J}}}
\newcommand{\KKK}{{\mathbf{K}}}
\newcommand{\LLL}{{\mathbf{L}}}
\newcommand{\MMM}{{\mathbf{M}}}
\newcommand{\NNN}{{\mathbf{N}}}
\newcommand{\OOO}{{\mathbf{O}}}
\newcommand{\PPP}{{\mathbf{P}}}
\newcommand{\QQQ}{{\mathbf{Q}}}
\newcommand{\RRR}{{\mathbf{R}}}
\newcommand{\SSS}{{\mathbf{S}}}
\newcommand{\TTT}{{\mathbf{T}}}
\newcommand{\UUU}{{\mathbf{U}}}
\newcommand{\VVV}{{\mathbf{V}}}
\newcommand{\WWW}{{\mathbf{W}}}
\newcommand{\XXX}{{\mathbf{X}}}
\newcommand{\YYY}{{\mathbf{Y}}}
\newcommand{\ZZZ}{{\mathbf{Z}}}

\newcommand{\AAAA}{{\mathscr{A}}}
\newcommand{\BBBB}{{\mathscr{B}} }
\newcommand{\CCCC}{{\mathscr{C}} }
\newcommand{\DDDD}{{\mathscr{D}} }
\newcommand{\EEEE}{{\mathscr{E}} }
\newcommand{\FFFF}{{\mathscr{F}} }
\newcommand{\GGGG}{{\mathscr{G}}}
\newcommand{\HHHH}{{\mathscr{H}}}
\newcommand{\IIII}{{\mathscr{I}}}
\newcommand{\JJJJ}{{\mathscr{J}}}
\newcommand{\KKKK}{{\mathscr{K}}}
\newcommand{\LLLL}{{\mathscr{L}}}
\newcommand{\MMMM}{{\mathscr{M}}}
\newcommand{\NNNN}{{\mathscr{N}}}
\newcommand{\OOOO}{{\mathscr{O}}}
\newcommand{\PPPP}{{\mathscr{P}}}
\newcommand{\QQQQ}{{\mathscr{Q}}}
\newcommand{\RRRR}{{\mathscr{R}}}
\newcommand{\SSSS}{{\mathscr{S}}}
\newcommand{\TTTT}{{\mathscr{T}}}
\newcommand{\UUUU}{{\mathscr{U}}}
\newcommand{\VVVV}{{\mathscr{V}}}
\newcommand{\WWWW}{{\mathscr{W}}}
\newcommand{\XXXX}{{\mathscr{X}}}
\newcommand{\YYYY}{{\mathscr{Y}}}
\newcommand{\ZZZZ}{{\mathscr{Z}}}

\newcommand{\gA}{{\mathfrak{A}}}
\newcommand{\gB}{{\mathfrak{B}}}
\newcommand{\gC}{{\mathfrak{C}}}
\newcommand{\gD}{{\mathfrak{D}}}
\newcommand{\gE}{{\mathfrak{E}}}
\newcommand{\gF}{{\mathfrak{F}}}
\newcommand{\gG}{{\mathfrak{G}}}
\newcommand{\gH}{{\mathfrak{H}}}
\newcommand{\gI}{{\mathfrak{I}}}
\newcommand{\gJ}{{\mathfrak{J}}}
\newcommand{\gK}{{\mathfrak{K}}}
\newcommand{\gL}{{\mathfrak{L}}}
\newcommand{\gM}{{\mathfrak{M}}}
\newcommand{\gN}{{\mathfrak{N}}}
\newcommand{\gO}{{\mathfrak{O}}}
\newcommand{\gP}{{\mathfrak{P}}}
\newcommand{\gQ}{{\mathfrak{Q}}}
\newcommand{\gR}{{\mathfrak{R}}}
\newcommand{\gS}{{\mathfrak{S}}}
\newcommand{\gT}{{\mathfrak{T}}}
\newcommand{\gU}{{\mathfrak{U}}}
\newcommand{\gV}{{\mathfrak{V}}}
\newcommand{\gW}{{\mathfrak{W}}}
\newcommand{\gX}{{\mathfrak{X}}}
\newcommand{\gY}{{\mathfrak{Y}}}
\newcommand{\gZ}{{\mathfrak{Z}}}

\newcommand{\gota}{{\mathfrak{a}}}
\newcommand{\gotb}{{\mathfrak{b}}}
\newcommand{\gotc}{{\mathfrak{c}}}
\newcommand{\gotd}{{\mathfrak{d}}}
\newcommand{\gote}{{\mathfrak{e}}}
\newcommand{\gotf}{{\mathfrak{f}}}
\newcommand{\gotg}{{\mathfrak{g}}}
\newcommand{\goth}{{\mathfrak{h}}}
\newcommand{\goti}{{\mathfrak{i}}}
\newcommand{\gotj}{{\mathfrak{j}}}
\newcommand{\gotk}{{\mathfrak{k}}}
\newcommand{\gotl}{{\mathfrak{l}}}
\newcommand{\gotm}{{\mathfrak{m}}}
\newcommand{\gotn}{{\mathfrak{n}}}
\newcommand{\goto}{{\mathfrak{o}}}
\newcommand{\gotp}{{\mathfrak{p}}}
\newcommand{\gotq}{{\mathfrak{q}}}
\newcommand{\gotr}{{\mathfrak{r}}}
\newcommand{\gots}{{\mathfrak{s}}}
\newcommand{\gott}{{\mathfrak{t}}}
\newcommand{\gotu}{{\mathfrak{u}}}
\newcommand{\gotv}{{\mathfrak{v}}}
\newcommand{\gotw}{{\mathfrak{w}}}
\newcommand{\gotx}{{\mathfrak{x}}}
\newcommand{\goty}{{\mathfrak{y}}}
\newcommand{\gotz}{{\mathfrak{z}}}



\newcommand{\kkrest}{\begin{picture}(14,14)
\put(00,04){\line(1,0){14}}
\put(00,02){\line(1,0){14}}
\put(06,-4){\line(0,1){14}}
\put(08,-4){\line(0,1){14}}
\end{picture}     }

\newcommand{\krest}{~\kkrest~}


\newcommand{\grd}{{\text{\rm grd}}}
\newcommand{\id}{\text{id}}
\newcommand{\Tb}{\text{ \rm Tb}}
\newcommand{\Log}{\text{\rm Log }}
\newcommand{\Wh}{\text{\rm Wh}}
\newcommand{\Ker}{\text{\rm Ker }}
\newcommand{\Ext}{\text{\rm Ext}}
\newcommand{\Hom}{\text{\rm Hom}}
\newcommand{\diam}{\text{\rm diam}}
\newcommand{\Homb}{\text{\rm Hom}b}
\newcommand{\Lg}{\text{\rm Lg }}
\newcommand{\ind}{\text{\rm ind\hspace{0.05cm}}}
\newcommand{\rk}{\text{\rm rk }}
\renewcommand{\Im}{\text{\rm Im }}
\newcommand{\supp}{\text{\rm supp }}
\newcommand{\Int}{\text{\rm Int }}
\newcommand{\grad}{\text{\rm grad}}
\newcommand{\Fix}{\text{\rm Fix}}
\newcommand{\Exp}{\text{\rm Exp}}
\newcommand{\Per}{\text{\rm Per}}
\newcommand{\TL}{\text{\rm TL}}
\newcommand{\Id}{\text{\rm Id}}
\newcommand{\Vect}{\text{\rm Vect}}
\newcommand{\vvol}{\text{\rm vol}}
\newcommand{\Mat}{\text\rm Mat}
\newcommand{\Tub}{\text{\rm Tub}}
\newcommand{\Imm}{\text{\rm Im}}
\newcommand{\tn}{\text{\rm t.n.}}
\newcommand{\card}{\text{\rm card }}
\newcommand{\GL}{\text{\rm GL}}

\newcommand{\track}{\text{\rm Track}}
\newcommand{\sgn}{\text{\rm sgn}}
\newcommand{\Arctg}{\text{\rm Arctg }}
\newcommand{\Arcsin}{\text{\rm Arcsin }}
\newcommand{\Det}{\text{\rm Det}}
\newcommand{\eddd}{\'equation diff\'erentielle}


\newcommand{\bere}{\begin{rema}}
\newcommand{\bede}{\begin{defi}}

\renewcommand{\beth}{\begin{theo}}
\newcommand{\bele}{\begin{lemm}}
\newcommand{\bepr}{\begin{prop}}
\newcommand{\beeq}{\begin{equation}}
\newcommand{\bega}{\begin{gather}}
\newcommand{\begaa}{\begin{gather*}}
\newcommand{\been}{\begin{enumerate}}

\newcommand{\bedee}{\begin{defii}}
\newcommand{\bethh}{\begin{theoo}}
\newcommand{\belee}{\begin{lemmm}}
\newcommand{\beprr}{\begin{propp}}

\newcommand{\beco}{\begin{coro}}

\newcommand{\beal}{\begin{aligned}}

\newcommand{\enre}{\end{rema}}

\newcommand{\enco}{\end{coro}}
\newcommand{\enpr}{\end{prop}}
\newcommand{\enth}{\end{theo}}
\newcommand{\enle}{\end{lemm}}
\newcommand{\enen}{\end{enumerate}}
\newcommand{\enga}{\end{gather}}
\newcommand{\engaa}{\end{gather*}}
\newcommand{\eneq}{\end{equation}}
\newcommand{\enal}{\end{aligned}}

\newcommand{\bq}{\begin{equation}}
\newcommand{\bqq}{\begin{equation*}}


\renewcommand{\leq}{\leqslant}
\renewcommand{\geq}{\geqslant}
\newcommand{\vphi}{\varphi}
\newcommand{\vide}{  \emptyset  }
\newcommand{\bu}{\bullet}
\newcommand{\pfff}{\pitchfork}
\newcommand{\mx}{\mbox}

\newcommand{\mxx}[1]{\quad\mbox{#1}\quad}
 \newcommand{\mxxx}[1]{\hspace{0.1cm}\mbox{#1} \quad  }
\newcommand{\wi}{\widetilde}

\newcommand{\ove}{\overline}
\newcommand{\unde}{\underline}
\newcommand{\ptf}{\pitchfork}

\newcommand{\emp}{\emptyset}
\newcommand{\wh}{\widehat}

\newcommand{\sub}{ Subsection}

\newcommand{\lc}{\lceil}
\newcommand{\rc}{\rceil}
\newcommand{\sps}{\supset}

\newcommand{\sm}{\setminus}
\newcommand{\ems}{\varnothing}
\newcommand{\sbs}{\subset}

\newcommand{\subs}{\subsection}
\newcommand{\ity}{\infty}


\newcommand{\GC}{\GG\CC}
\newcommand{\GCT}{\GG\CC\TT}
\newcommand{\GT}{\GG\TT}

\newcommand{\GA}{\GG\AA}
\newcommand{\GRP}{\GG\RR\PP}

\newcommand{\GgC}{\GG\gC}
\newcommand{\GgCC}{\GG\gC\CC}

\newcommand{\GgCT}{\GG\gC\TT}

\newcommand{\GgCY}{\GG\gC\YY}
\newcommand{\GgCYT}{\GG\gC\YY\TT}

\newcommand{\GCCT}{\GG\gC\CC\TT}
\newcommand{\GCC}{\GG\gC\CC}

\newcommand{\GKSC}{\GG\KK\SS\gC}
\newcommand{\GKS}{\GG\KK\SS}


\newcommand{\dstr}[1]
{
\TT_{#1}
}

\newcommand{\strr}[3]
{{#1}^{\displaystyle\twoheadrightarrow}_{[{#2},{#3}]}}

\newcommand{\str}[1]{{#1}^{\displaystyle\twoheadrightarrow}}

\newcommand{\stind}[3]
{{#1}^{\displaystyle\rightsquigarrow}_{[{#2},{#3}]}}

\newcommand{\st}[1]{\overset{\rightsquigarrow}{#1}}
\newcommand{\bst}[1]{\overset{\displaystyle\rightsquigarrow}
\to{\boldkey{#1}}}

\newcommand{\stexp}[1]{{#1}^{\rightsquigarrow}}
\newcommand{\bstexp}[1]{{#1}^{\displaystyle\rightsquigarrow}}

\newcommand{\bstind}[3]{{\boldkey{#1}}^{\displaystyle\rightsquigarrow}_
{[{#2},{#3}]}}
\newcommand{\bminstind}[3]{\stind{({\boldkey{-}\boldkey{#1}})}{#2}{#3}}

\newcommand{\ST}{\stexp}

\newcommand{\stv}{\stexp {(-v)}}
\newcommand{\stu}{\stexp {(-u)}}
\newcommand{\stw}{\stexp {(-w)}}

\newcommand{\strv}[2]{\stind {(-v)}{#1}{#2}}
\newcommand{\strw}[2]{\stind {(-w)}{#1}{#2}}
\newcommand{\stru}[2]{\stind {(-u)}{#1}{#2}}

\newcommand{\stvv}[2]{\stind {v}{#1}{#2}}
\newcommand{\stuu}[2]{\stind {u}{#1}{#2}}
\newcommand{\stww}[2]{\stind {w}{#1}{#2}}

\newcommand{\vovo}{\stexp {(-v0)}}
\newcommand{\vov}{\stexp {(-v1)}}

\newcommand{\fl}[1]{{#1}\!\da}
\newcommand{\fll}[1]{({#1}\!\da)}

\newcommand{\vflesh}{\fl{v}}
\newcommand{\wflesh}{\fl{w}}

\newcommand{\vfllesh}{\fll{v}}
\newcommand{\wfllesh}{\fll{w}}

\newcommand{\vvfl}{\wi v\!\da}
\newcommand{\wivflesh}{\wi v\!\da}

\newcommand{\RA}{\Rightarrow}
\newcommand{\LA}{\Leftarrow}
\newcommand{\RLA}{\Leftrightarrow}

\newcommand{\LRA}{\Leftrightarrow}

\newcommand{\lau}[1]{{\xleftarrow{#1}}}

\newcommand{\rau}[1]{{\xrightarrow{#1}}}
\newcommand{\rad}[1]{ {\xrightarrow[#1]{}} }

\newcommand{\da}{\downarrow}


\newcommand{\vecm}{\Vect^1_0(M)}
\newcommand{\vecbn}{\Vect^1_b(N)}
\newcommand{\vecw}{\Vect^1(W)}
\newcommand{\vecrm}{\Vect^1_b(\RRR^m)}
\newcommand{\vecbm}{\Vect^1_b(M)}
\newcommand{\ver}{\text{\rm Vect}^1(\RRR^ n)}
\newcommand{\verr}{\text{\rm Vect}^1_0(\RRR^ n)}
\newcommand{\hrrr}{\text{\rm Vect}^1(M)}
\newcommand{\vemm}{\text{\rm Vect}^1_0(M)}

\newcommand{\vem}{\text{{\rm Vectt}}(M)}
\newcommand{\vebK}{\text{{\rm Vectt}}(B,K)}
\newcommand{\vemK}{\text{{\rm Vectt}}(M,K)}
\newcommand{\vemc}{\text{{\rm Vectt}}_c(M)}
\newcommand{\vemQ}{\text{{\rm Vectt}}(M,Q)}

\newcommand{\vectt}[1]{\text{{\rm Vectt}}(#1)}

\newcommand{\vecsmo}[1]{{\text{\rm Vect}}^\infty (#1)}

\newcommand{\vew}{\text{\rm Vect}^1 (W,\bot)}

\newcommand{\downnorm}{\text{\rm Vect}^1_\bot (W)}

\newcommand{\upnorm}{\text{\rm Vect}^1_\top (W)}

\newcommand{\normm }{\text{\rm Vect}^1_N (W)}

\newcommand{\veww}{\text{\rm Vect}^1 (W)}


\newcommand{\tens}[1]{\underset{#1}{\otimes}}

\newcommand{\starr}[1]{\underset{#1}{*}}

\newcommand{\Lxxxi}{\wh L_{\bar\xi}}

\newcommand{\Lbarxi}{\wh \L_{\bar\xi}}

\newcommand{\RRxi}{\wh \RR_\xi}
\newcommand{\LLxi}{\wi \L_\xi}

\newcommand{\Lx}{\L_{(\xi)}}
\newcommand{\Lxi}{{\wh \L}_\xi}
\newcommand{\Leta}{{\wh \L}_\eta}

\newcommand{\lL}{\wh{\wh L}}

\newcommand{\Rxi}{{\ove R}_\xi}
\newcommand{\Nxi}{{\ove N}_\xi}
\newcommand{\Rcxi}{{\bar R}_\xi^c}

\newcommand{\sil}{ S^{-1}\L }
\newcommand{\kil}{\ove{K}_1(\L)}
\newcommand{\killl}{\ove{K}_1(\wh\L)}
\newcommand{\kisl}{\ove{K}_1(S^{-1}\L )}

\newcommand{\klxi}{K_1(\Lxi)}
\newcommand{\kklxi}{\ove{K_1}(\Lxi)}

\newcommand{\popo}{\tens{\L}\Lxi}
\newcommand{\popom}{\tens{\L^-}\Lxi^-}

\newcommand{\wll}{{\text{\rm Wh}}\big(\hspace{0.06cm} \ove L\hspace{0.06cm}\big)}
\newcommand{\kll}{ K_1\big(\hspace{0.06cm} \ove L\hspace{0.06cm}\big)}


\newcommand{\amk}{\AA^{(m)}_k}
\newcommand{\amkm}{\AA^{(m)}_{k-1}}
\newcommand{\bmk}{\BB^{(m)}_k}
\newcommand{\bmkm}{\BB^{(m)}_k}

\newcommand{\tivm}{\wi V^-}

\newcommand{\vk}{V_{\langle k\rangle}^-}
\newcommand{\tivkm}{\wi V_{\langle k-1\rangle}^-}
\newcommand{\tivk}{\wi V_{\langle k\rangle}^-}

\newcommand{\vkm}{V_{\langle k-1\rangle}^-}
\newcommand{\vkp}{V_{\langle k+1\rangle}^-}

\newcommand{\hvk}{\wh V_{\langle k\rangle}^-}

\newcommand{\hvkm}{\wh V_{\langle k-1\rangle}^-}
\newcommand{\hvkp}{\wh V_{\langle k+1\rangle}^-}

\newcommand{\vkvk}{V_{\prec k\succ}^-}

\newcommand{\vkvkm}{V_{\prec k-1\\succ}^-}

\newcommand{\tivkvk}{\wi V_{\prec k\succ}^-}

\newcommand{\tivkvkm}{\wi V_{\prec k-1 \succ}^-}

\newcommand{\vvvbs}{V_b^{( s)    }}
\newcommand{\vvvas}{V_a^{( s)    }}
\newcommand{\vvvbsm}{V_b^{( s-1)    }}
\newcommand{\vvvasm}{V_a^{( s-1)    }}

\newcommand{\vvbsm}{V_b^{( s-1)    }}
\newcommand{\vvasm}{V_a^{( s-1)    }}
\newcommand{\vvbs}{V_b^{[s]}    }
\newcommand{\vvas}{V_a^{[ s]    }}

\newcommand{\factor}{\vvbs / \vvbsm}
\newcommand{\factora}{\vvas / \vvasm}

\newcommand{\vvksm}{\wi V_k^{( s-1)   }}
\newcommand{\vvks}{\wi V_k^{[s]}    }

\newcommand{\vvkmsm}{\wi V_{k-1}^{( s-1)   }}
\newcommand{\vvkms}{\wi V_{k-1}^{[s]}    }

\newcommand{\fac}{\vvks / \vvksm}
\newcommand{\facm}{\vvkms / \vvkmsm}

\newcommand{\vbsm}{V_b^{\{\leq s-1\}}    }
\newcommand{\vasm}{V_a^{\{\leq s-1\}}    }
\newcommand{\vbs}{V_b^{\{\leq s\}}    }
\newcommand{\vas}{V_a^{\{\leq s\}}    }

\newcommand{\wivksm}{\wi V_k^{\{\leq s-1\}}    }
\newcommand{\wivkmsm}{\wi V_{k-1}^{\{\leq s-1\}}    }
\newcommand{\wivks}{\wi V_k^{\{\leq s\}}    }
\newcommand{\wivkms}{\wi V_{k-1}^{\{\leq s\}}    }

\newcommand{\Vbsm}{V_b^{[\leq s-1]}(\d)    }
\newcommand{\Vasm}{V_a^{[\leq s-1]}(\d)    }
\newcommand{\Vbs}{V_b^{[\leq s]}(\d)    }
\newcommand{\Vas}{V_a^{[\leq s]}(\d)    }

\newcommand{\vass}{V_{a_{s+1}}}

\newcommand{\vbkm}{V_b^{\{\leq k-1\}}    }
\newcommand{\vakm}{V_a^{\{\leq k-1\}}    }
\newcommand{\vbk}{V_b^{\{\leq k\}}    }
\newcommand{\vak}{V_a^{\{\leq k\}}    }

\newcommand{\Vbkm}{V_b^{[\leq k-1]}(\d)    }
\newcommand{\Vakm}{V_a^{[\leq k-1]}(\d)    }
\newcommand{\Vbk}{V_b^{[\leq k]}(\d)    }
\newcommand{\Vak}{V_a^{[\leq s]}(\d)    }


\newcommand{\dow}{\pr_0 W}

\newcommand{\daw}{\pr_1 W}

\newcommand{\hdaw}{\wh{\pr_1 W}}

\newcommand{\tipwk}{(\pr_1 \wi W)^{\{\leq k\}}}

\newcommand{\timwk}{(\pr_0 \wi W)^{\{\leq k\}}}

\newcommand{\tipwkm}{(\pr_1 \wi W)^{\{\leq k-1\}}}

\newcommand{\pws}{(\pr_1 W)^{\{\leq s\}}}

\newcommand{\hpws}{(\wh{\pr_1 W})^{\{\leq s\}}}

\newcommand{\hpwsm}{(\wh{\pr_1 W})^{\{\leq {s-1}\}}}

\newcommand{\mws}{(\pr_0 W)^{\{\leq s\}}}

\newcommand{\pwk}{(\pr_1 W)^{\{\leq k\}}}

\newcommand{\mwk}{(\pr_0 W)^{\{\leq k\}}}

\newcommand{\pwkm}{(\pr_1 W)^{\{\leq k-1\}}}

\newcommand{\pwsm}{(\pr_1 W)^{\{\leq s-1\}}}

\newcommand{\pwkmm}{(\pr_1 W)^{\{\leq k-2\}}}

\newcommand{\mwkmm}{(\pr_0 W)^{\{\leq k-2\}}}

\newcommand{\mwkm}{(\pr_0 W)^{\{\leq k-1\}}}

\newcommand{\mwsm}{(\pr_0 W)^{\{\leq s-1\}}}

\newcommand{\mwkp}{(\pr_0 W)^{\{\leq k+1\}}}

\newcommand{\dwmok}{\daw^{( k)}}

\newcommand{\dwmokp}{\daw^{(k+1)}}

\newcommand{\dwmokm}{\daw^{(k-1)}}
\newcommand{\dwmokmm}{\daw^{(k-2)}}


\newcommand{\moi}[1]{\MM^{(0)}_{#1}}
\newcommand{\moii}[1]{\MM^{(1)}_{#1}}

\newcommand{\Wal}{W_{[a,\l]}}

\newcommand{\Wlm}{W_{[\l,\m]}}
\newcommand{\Wlb}{W_{[\l,b]}}

\newcommand{\Wam}{W_{[a,\m]}}

\newcommand{\Wall}{W_{[a,\l']}}

\newcommand{\Wkr}{W^{\circ}}
\newcommand{\wkr}{W^{\circ}}

\newcommand{\wa}[2]{ W_{[a_{#1}, a_{#2}]}}

\newcommand{\waa}[1]{ W_{[a, a_{#1}]}}

\newcommand{\Wa}[2]{ W_{[{#1}, {#2}]}}

\newcommand{\WS}[1]{ W^{\{\leq {#1}\}}}

\newcommand{\ws}{\WS {s}}

\newcommand{\wsm}{\WS {s-1}}

\newcommand{\wsmm}{\WS {s-2}}

\newcommand{\wk}{\WS {k}}

\newcommand{\wkm}{\WS {k-1}}

\newcommand{\wkmm}{\WS {k-2}}

\newcommand{\wsn}{ W^{[\leq s]}(\nu)}

\newcommand{\wsmn}{ W^{[\leq s-1]}(\nu)}

\newcommand{\wsk}{ W^{[\leq k]}(\nu)}

\newcommand{\Wmok}{W^{\prec k\succ}}
\newcommand{\wmok}{W^{\langle k\rangle}}

\newcommand{\wmokp}{W^{\langle k+1\rangle}}

\newcommand{\wmokm}{W^{\langle k-1\rangle}}
\newcommand{\wmokmm}{W^{\langle k-2\rangle}}

\newcommand{\wmos}{W^{\langle s\rangle}}

\newcommand{\wmosm}{W^{\langle s-1\rangle}}

\newcommand{\wmoo}{W^{\langle 0\rangle}}

\newcommand{\wwk}{\( \wmok , \wmokm \)}

\newcommand{\wwkp}{\( \wmokp , \wmok \)}

\newcommand{\wwkm}{\( \wmokm , \wmokmm \)}

\newcommand{\wws}{\bigg( \wmos , \wmosm \bigg)}

\newcommand{\wasn}{W^{\lc s\rc}(\nu)}

\newcommand{\wakn}{W^{\lc k\rc}(\nu)}

\newcommand{\hwm}{H_*\( \wmok, \wmokm\)}

\newcommand{\hkwm}{H_k\( \wmok, \wmokm\)}

\newcommand{\hwmp}{H_*\( \wmokp, \wmok\)}

\newcommand{\hkwmm}{H_k\( \wmokm, \wmokmm\)}

\newcommand{\hkwmp}{H_k\( \wmokp, \wmok\)}

\newcommand{\tiws}{\wi W^{\{\leq s\}}}

\newcommand{\tiwk}{\wi W^{\{\leq k\}}}

\newcommand{\tiwsm}{\wi W^{\{\leq s-1\}}}

\newcommand{\tiwkm}{\wi W^{\{\leq k-1\}}}

\newcommand{\wkwk}{ W^{[ k]}}
\newcommand{\tiwmok}{\wi W^{\prec k\succ}}
\newcommand{\womok}{W_0^{\prec k\succ}}

\newcommand{\Womok}{W_0^{\langle k\rangle}}

\newcommand{\Wmokp}{W^{\langle k+1\rangle}}

\newcommand{\Wmokm}{W^{\langle k-1\rangle}}
\newcommand{\Wmokmm}{W^{\langle k-2\rangle}}

\newcommand{\Wmos}{W^{\langle s\rangle}}

\newcommand{\Wmosm}{W^{\langle s-1\rangle}}

\newcommand{\Wmoo}{W^{\langle 0\rangle}}

\newcommand{\hWmok}{\wh W^{\langle k\rangle}}
\newcommand{\hWmokm}{\wh W^{\langle k-1\rangle}}

\newcommand{\tiWmok}{\wi W^{\langle k\rangle}}
\newcommand{\tiWmokm}{\wi W^{\langle k-1\rangle}}



\newcommand{\talp}{twisted Alexander polynomial}
\newcommand{\tnh}{twisted Novikov homology}

\newcommand{\ifff}{if and only if}

\newcommand{\orial}{oriented almost transverse}
\renewcommand{\th}{therefore}
\newcommand{\at}{almost~ transverse}
\newcommand{\ata}{almost~ transversality~ condition}
\newcommand{\gr}{gradient}
\newcommand{\Mf}{Morse function}
\newcommand{\iis}{it is sufficient}
\newcommand{\sut}{~such~that~}
\newcommand{\sufsm}{~sufficiently~ small}
\newcommand{\sufla}{~sufficiently~ large}
\newcommand{\sufcl}{~sufficiently~ close}
\newcommand{\wrt}{with respect to}
\newcommand{\ho}{homomorphism}
\newcommand{\iso}{isomorphism}
\newcommand{\rgr}{Riemannian gradient}
\newcommand{\riemm}{Riemannian metric}

\newcommand{\trasp}{trajectory starting at a point of}
\newcommand{\trasps}{trajectories starting at a point of}

\newcommand{\ma}{manifold}
\newcommand{\nei}{neighbourhood}
\newcommand{\dfm}{diffeomorphism}

\newcommand{\vf}{vector field}

\newcommand{\vfs}{vector fields}

\newcommand{\fe}{for every}

\newcommand{\tr}{~trajectory }

\newcommand{\grs}{~gradients}
\newcommand{\trs}{~trajectories}

\newcommand{ \co}{~cobordism}
\newcommand{
\sma}{submanifold}
\newcommand{
\hos}{~homomorphisms}
\newcommand{
\Th}{~Therefore}

\newcommand{
\tthen}{\text \rm ~then}

\newcommand{
\wwe}{\text \rm ~we  }
\newcommand{
\hhave}{\text \rm ~have}
\newcommand{
\eevery}{\text \rm ~every}

\newcommand{\noconf}{~no~confusion~is~possible}

\newcommand{\ATA}{almost~ transversality~ condition}
\newcommand{\cob}{~cobordism}

\newcommand{\hot}{homotopy}

\newcommand{\emem}{elementary modification}
\newcommand{\emems}{elementary modifications}

\newcommand{\TA}{transversality condition}

\newcommand{\hog}{homology}

\newcommand{\cog}{cohomology}

\newcommand{\wat}{ We shall assume that}

\newcommand{\sclv}{sufficiently close to $v$ in $C^0$-topology}

\newcommand{\cf}{continuous function }

\newcommand{\heq}{homotopy equivalence}

\newcommand{\heeq}{homology equivalence}

\newcommand{\eg}{exponential growth}

\newcommand{\nics}{Novikov incidence coefficients}
\newcommand{\nic}{Novikov incidence coefficient}

\newcommand{\negc}{Novikov exponential growth conjecture}

\newcommand{\mas}{manifolds   }

\newcommand{\nc}{Novikov Complex   }

\newcommand{\glvf}{gradient-like vector field}

\newcommand{\glvfs}{gradient-like vector fields}

\newcommand{\fg}{finitely generated   }

\newcommand{\she}{simple~homotopy~equivalence}

\newcommand{\sht}{simple~homotopy~type}

\newcommand{\ta}{transversality condition}

\newcommand{\cpc}{convex polyhedral cone}
\newcommand{\rcpc}{rational convex polyhedral cone}

\newcommand{\mnp}{Morse-Novikov pair}

\newcommand{\rp}{rationality property}
\newcommand{\wvf}{Whitney vector field}

\newcommand{\egp}{exponential growth property}

\newcommand{\lzf}{Lefschetz zeta function}

\newcommand{\babs}{by abuse of notation}
\newcommand{\su}{subsection}
\newcommand{\Prop}{\text{Proposition}}

\newcommand{\aand}{\quad\text{and}\quad}
\newcommand{\wwhere}{\quad\text{where}\quad}
\newcommand{\ffor}{\quad\text{for}\quad}
\newcommand{\iif}{\quad\text{if}\quad}
\newcommand{\iiif}{~\text{if}~}

\newcommand{\eqi}{equivalence}

\newcommand{\mfcobv}{~ Let $\fcob$
be a Morse function on
 a cobordism $W$ and $v$
 an $f$-gradient. ~}

\newcommand{\mfcob}{~ Let $\fcob$
be a Morse function on
 a cobordism $W$}

\newcommand{\mfcobvat}{~ Let $\fcob$
be a Morse function on
 a cobordism $W$ and $v$
  an almost transverse $f$-gradient. ~}

\newcommand{\msf}{Morse-Smale filtration}

\newcommand{\fbfg}{ free based finitely generated }

\newcommand{\tap}
{twisted Alexander invariant}


\newcommand{\vaa}{\mathscr A_k}

\newcommand{\qaa}{\mathscr Q_k}

\newcommand{\tret}{{\frac 13}}
\newcommand{\dvet}{{\frac 23}}
\newcommand{\polt}{{\frac 32}}
\newcommand{\polo}{{\frac 12}}

\newcommand{\bv}{B(v,\d)}

\newcommand{\ti}{\times}

\newcommand{\FR}{{\mathcal{F}}r}

\newcommand{\en}{enumerate}

\newcommand{\Prf}{{\it Proof.\quad}}
\newcommand{\prf}{{\it Proof:\quad}}

\newcommand{\nr}{\Vert}
\newcommand{\smo}{C^{\infty}}

\newcommand{\fpr}[2]{{#1}^{-1}({#2})}
\newcommand{\sdvg}[3]{\widehat{#1}_{[{#2},{#3}]}}
\newcommand{\disc}[3]{B^{({#1})}_{#2}({#3})}
\newcommand{\Disc}[3]{D^{({#1})}_{#2}({#3})}
\newcommand{\desc}[3]{B_{#1}(\leq{#2},{#3})}
\newcommand{\Desc}[3]{D_{#1}(\leq{#2},{#3})}
\newcommand{\komp}[3]{{\bold K}({#1})^{({#2})}({#3})}
\newcommand{\Komp}[3]{\big({\bold K}({#1})\big)^{({#2})}({#3})}

\newcommand{\ran}{\{(A_\lambda , B_\lambda)\}_{\lambda\in\Lambda}}
\newcommand{\rran}{\{(A_\lambda^{(s)},
 B_\lambda^{(s)}  )\}_{\lambda\in\Lambda, 0\leq s\leq n }}
\newcommand{\brs}{\rran}
\newcommand{\rans}{\{(A_\sigma , B_\sigma)\}_{\sigma\in\Sigma}}

\newcommand{\fmin}{F^{-1}}
\newcommand{\vh}{\widehat{(-v)}}

\newcommand{\chart}{\Phi_p:U_p\to B^n(0,r_p)}
\newcommand{\atlas}{\{\Phi_p:U_p\to B^n(0,r_p)\}_{p\in S(f)}}
\newcommand{\flow}{{\VV}=(f,v, \UU)}

\newcommand{\Rn}{\bold R^n}
\newcommand{\Rk}{\bold R^k}

\newcommand{\fcob}{f:W\to[a,b]}

\newcommand{\phicob}{\phi:W\to[a,b]}

\newcommand{\crr}{p\in S(f)}
\newcommand{\nrv}{\Vert v \Vert}
\newcommand{\nrw}{\Vert w \Vert}
\newcommand{\nru}{\Vert u \Vert}

\newcommand{\obb}{\cup_{p\in S(f)} U_p}
\newcommand{\proob}{\Phi_p^{-1}(B^n(0,}

\newcommand{\indl}[1]{{\scriptstyle{\text{\rm ind}\leqslant {#1}~}}}
\newcommand{\inde}[1]{{\scriptstyle{\text{\rm ind}      =   {#1}~}}}
\newcommand{\indg}[1]{{\scriptstyle{\text{\rm ind}\geqslant {#1}~}}}

\newcommand{\obbi}{\cup_{p\in S_i(f)}}

\newcommand{\pr}{\partial}
\newcommand{\prx}[2]{\frac {\pr {#1}}{\pr x} ({#2})}
\newcommand{\pry}[2]{\frac {\pr {#1}}{\pr y} ({#2})}
\newcommand{\prz}[2]{\frac {\pr {#1}}{\pr z} ({#2})}
\newcommand{\przbar}[2]{\frac {\pr {#1}}{\pr \bar z} ({#2})}
\newcommand{\chape}[2]{\frac  {\pr {#1}}{\pr {#2}} }
\newcommand{\chapee}[2]{\frac  {\pr^2 {#1}}{\pr {#2}^2} }

\newcommand{\xit}{\tilde\xi_t}

\newcommand{\VODIN}{V_{1/3}}
\newcommand{\VDVA}{V_{2/3}}
\newcommand{\VM}{V_{1/2}}
\newcommand{\ddd}{\cup_{p\in S_i(F_1)} D_p(u)}
\newcommand{\dddmin}{\cup_{p\in S_i(F_1)} D_p(-u)}
\newcommand{\where}{\quad\text{\rm where}\quad}

\newcommand{\kr}[1]{{#1}^{\circ}}

\newcommand{\mods}{\vert s(t)\vert}
\newcommand{\exd}{e^{2(D+\alpha)t}}
\newcommand{\exmin}{e^{-2(D+\alpha)t}}

\newcommand{\intt}{[-\theta,\theta]}

\newcommand{\ffmin}{f^{-1}}

\newcommand{\vxi}{v\langle\vec\xi\rangle}

\newcommand{\qt}{\hfill\triangle}
\newcommand{\qs}{\hfill\square}

\newcommand{\pa}{\vskip0.1in}
\newcommand{\spp}{\hspace{0.07cm}}
\newcommand{\sppp}{\hspace{0.05cm}}
\newcommand{\spppp}{\hspace{0.03cm}}

\renewcommand{\(}{\big(}
\renewcommand{\)}{\big)}

\newcommand{\Vm}{V_\m}

\newcommand{\Vl}{V_\l}

\newcommand{\lccc}{\wh\L_{C}}

\newcommand{\ld}{\wh\L_{D}}

\newcommand{\udp}{{\displaystyle {\vartriangle}}}
\newcommand{\ddp}{{\displaystyle {\triangledown}}}

\newcommand{\Vv}{{\boldsymbol{v}}}

\newcommand{\hV}{\wh V}
\newcommand{\hHH}{\wh \HH}

\newcommand{\gama}[2]{\g({#1}, \tau_a({#2},{#1}); w )}

\newcommand{\gam}[2]{\g({#1}, \tau_0({#2},{#1}); w )}
\newcommand{\ga}[2]{\g({#1}, \tau({#2},{#1}); w )}

\newcommand{\mi}[3]{{#1}^{-1}\([{#2},{#3}]\)}

\newcommand{\fii}[2]{\mi {\phi}{a_{#1}}{a_{#2}} }

\newcommand{\fifi}[2]{\mi {\phi}{#1}{#2} }

\newcommand{\mf}[2]{\mi {\phi_0}{\b_{#1}}{\b_{#2}}}

\newcommand{\dqr}{\pr_- Q_r}

\newcommand{\ds}{\pr_s}

\newcommand{\dsm}{\pr_{s-1}}

\newcommand{\yz}{Y_k(v)\cup Z_k(v)}

\newcommand{\Gama}{{\nazad{ \Gamma}}}
\newcommand{\ug}[1]{\llcorner {#1} \lrcorner}
\newcommand{\npqv}{n(\bar p, \bar q; v)}
\newcommand{\fms}{f:M\to S^1   }

\newcommand{\nkpqv}{n_k(\bar p, \bar q; v)}

\newcommand{\GLT}{\GG lt}

\newcommand{\Trln}{{\text Trln}}

\newcommand{\Trlln}{{\text TrLn}}

\newcommand{\Tr}{{\text{\rm  Tr}}}
\newcommand{\TrL}{{\text TrL}}
\newcommand{\limdir}{\underset {\to}{\lim}}
\newcommand{\liminv}{\underset {\leftarrow}{\lim}}

\newcommand{\kom}[2]{ {#1}{#2}{ {#1}^{-1}} {{#2}^{-1}} }

\newcommand{\komm}[2]{ {#1}{#2}{ ({#1})^{-1}} {({#2})^{-1}} }
\newcommand{\kommm}[2]{ {#1}'{#2}'{ ({#1}'')^{-1}} {({#2}'')^{-1}} }

\newcommand{\Trll}{\TL'}

\newcommand{\cmd}{ C_*^\D( \wi M)}
\newcommand{\cmxi}{\wh C_*^\D( \wi M, \xi)}
\newcommand{\whgxi}{\wh {{\rm Wh}} (G,\xi)}

\newcommand{\vwdwp}{Vect(W,\pr_0W;P)}
\newcommand{\ewdwp}{\EE(W,\pr_0W)}
\newcommand{\ewdwo}{\EE(W_1,\pr_0W_1)}
\newcommand{\ewdwd}{\EE(W_2,\pr_0W_2)}

\newcommand{\kpr}{K_r^+}

\newcommand{\kmr}{K_r^-}

\newcommand{\kpd}{K_r^+(\d)}

\newcommand{\kmd}{K_r^-(\d)}

\newcommand{\addc}{\addtocontents{toc}{\protect\vspace{10pt}}}

\newcommand{\mxi}{M_\xi   }

\newcommand{\cmm}{C_*^\D(\wi M)}

\newcommand{\cvm}{C_*^\D(\wi V^-)}

\newcommand{\ey}{\wi E_*}
\newcommand{\eey}{\wi \EE_*}

\newcommand{\eky}{\wi E(k)_*}

\newcommand{\eti}{\wi{\wi\EE_*}}

\newcommand{\etik}{\wi{\wi\EE}_k}

\newcommand{\etikp}{\wi{\wi\EE}_{k+1}}

\newcommand{\ctiu}{\wi{\wi C}_*(u_1)}

\newcommand{\ctiuk}{\wi{\wi C}_k(u_1)}

\newcommand{\ctiv}{\wi{\wi C}_*(v)}

\newcommand{\ctiukm}{\wi{\wi C}_{k-1}(u_1)}

\newcommand{\scc}[1]{|{\scriptscriptstyle{#1}}}\newcommand{\rrr}{\{\wi r\}}

\newcommand{\tidow}{\pr_0 \wi W}
\newcommand{\tidaw}{\pr_1 \wi W}

\newcommand{\tivkp}{\wi V_{\langle k+1\rangle}^-}

\newcommand{\ur}[1]{\overset{\smallfrown}{#1}}

\newcommand{\dr}[1]{\underset{\smallsmile}{#1}}

\newcommand{\uUu}{\overset{\twoheadrightarrow}{u}}
\newcommand{\vVv}{\overset{\twoheadrightarrow}{v}}
\newcommand{\wWw}{\overset{\twoheadrightarrow}{w}}

\newcommand{\bfun}{{\bf 1}}

\newcommand{\ppmm}{{\scriptstyle{ \pm}}}

\newcommand{\dpm}{{\scriptstyle{ \pm}}}


\newcommand{\Lxim}{{\wh \L}^-_\xi}
\newcommand{\Lgxi}{{\wh \gL}_\xi}
\newcommand{\lxi}{{\bar \L}_\xi}

\newcommand{\Lc}{{\wh \L}_C}
\newcommand{\Lcm}{{\wh \L}_C^-}

\newcommand{\lLL}{\wh{\wh \L}}

\newcommand{\Xc}{{\wh X}_C}
\newcommand{\Xfaa}{{\wh X}_{(F_i,\vec\a)}}

\newcommand{\bs}{\boldsymbol}


\newcommand{\bikl}{\text{\rM (}}
\newcommand{\bikr}{\text{\rM )}}

\newcommand{\ck}{C_*^{(k)}}
\newcommand{\dk}{D_*^{(k)}}
\newcommand{\ckp}{C_*^{(k+1)}}
\newcommand{\dkp}{D_*^{(k+1)}}

\newcommand{\tlft}{\vartriangleleft}

\newcommand{\AbrahamRobbinTransv}
{   R. Abraham, J.Robbin,
\emph{Transversal mappings and flows
},
 Benjamin, New York,  1967.
}

\newcommand{\AlaniyaAl}
{L.Alaniya,
\emph{ Manifolds of Alexander type} (in Russian),
Uspekhi Mat.Nauk, {\bff 46 } N.1, 1991.
\quad English Transl:
Russian Math.Surveys, {\bff 46 } N.1, 1991.}

\newcommand{\ArnoldDynInt}
{V.I.Arnold,
\emph{ Dynamics of intersections},
 Proceedings of a Conference
 in Honour of
J.Moser, edited by
 P.Rabinowitz and R.Hill,
 Academic Press,
 1990
pp. 77--84.
}



\newcommand{\ArnoldEquadiff}
{   V.I.Arnold,
\emph{Ordinary Differential Equations},
 Moscow, Nauka,  1975.
}

\newcommand{\ArnoldDynComplInt}
{ V.I.Arnold,
\emph{Dynamics of Complexity of Intersections},
Boletim Soc. Brasil. Mat. (N.S.),
1990, {\bff 21} (1), 1-10.
}



\newcommand{\ArnoldBoundsMiln}{ 
V.I. Arnold,
\emph{Bounds for Milnor numbers of intersections
in holomorphic dynamical systems}, In:
Topological Methods in Modern Mathematics,
 Publish or Perish,
 1993,
pp. 379--390.}



\newcommand{\ArnoldProbl}{ V.I. Arnold,
\emph{Sur quelques probl\`emes de
 la th\'eorie des syst\`emes
dynamiques},
  Journal of the Julius Schauder center, {\bff 4}
 1994
pp. 209--225.  }



\newcommand{\ArtinMazurPerPo}
{M.Artin, B.Mazur,
\emph{On periodic points},
Annals of Math. {\bff 102} (1965), 82--99.
}



\newcommand{\AtiyahMacdonaldCA}
{   M.F.Atiyah, I.G.Macdonald, 
\emph{Introduction to commutative algebra},
Addison-Wesley,   1969.
}

\newcommand{\BaladiPerOrb}
{V.Baladi,
\emph{Periodic orbits and dynamical spectra},
Ergodic theory and dynamical systems,
{\bff 18}
(1998),
    255 - 292.
}



\newcommand{\BassAlgK}
{ H.Bass,
\emph{Algebraic K-theory},
Benjamin, 1968.
}

\newcommand{\BassHellerSwanWh}
{H.Bass, A.Heller, R.G.Swan,
\emph{The Whitehead group of a polynomial extension},
Inst. Hautes Etudes Sci. Publ. Math. {\bff 22} (1964),
61--79
}

\newcommand{\BirkhoffRotaDiffEq}
{ G.Birkhoff, G-C. Rota,
\emph{Ordinary differential equations},
Blaisdell Publishing Company, 1962.
}

\newcommand{\BourbakiGrLie}
{N.Bourbaki,
\emph{Groupes de Lie, Alg\`ebres de Lie},
Hermann, 1971.
}

\newcommand{\BravermanFarberNTineq}
{M. Braverman, M. Farber,
\emph{ Novikov type inequalities for differential forms with non-isolated zeros},
e-print: dg-ga/9508006, Journal publication:
 Math. Proc. Cambridge Philos. 
 Soc., 122 (1997) 357-375.
}

\newcommand{\BravermanFarberNBineq}
{M. Braverman, M. Farber,
\emph{ The Novikov-Bott inequalities},
 C.R. Acad. Sci. Paris, t. 321 (1995), 897-902.
}

\newcommand{\BrowderLevineFiber}
{ W.Browder, J.Levine,
\emph{Fibering manifolds over the circle},
Comment. Math. Helv. {\bff 40} (1966),
153--160.
}

\newcommand{\BurdeZieschang}
{ G.Burde, H.Zieschang,
\emph{ Knots
}, Walter de Gruyter,  1985.}

\newcommand{\BurgheleaHallerNovRev}
{ D.Burghelea, S.Haller, 
\emph{ On the topology and analysis of a closed 1-form
(Novikov's theory revisited)
}, \quad  E-print: dg-ga/0101043
Journal publication: 
Essays on geometry and related topics, Vol.1,2, 133-175,
Monogr. Enseign. Math., 38,
Enseignement Math., Geneva, 2001.}

\newcommand{\CartanCD}
{   H.Cartan,
\emph{Cours de Calcul Diff\'erentiel},
Hermann,  1977.
}

\newcommand{\ChapmanTopInv}
{ T.A.Chapman,
\emph{ Topological invariance of Whitehead torsion},
Amer J. Math.
{\bff 96}
(1974),
    488 - 497.
}
\newcommand{\CohenSHT}
{
M.M.Cohen,
\emph{A course in Simple-Homotopy theory},
Springer, 1972.}



\newcommand{\CohnFreeRings}
{ P.M.Cohn,
\emph{Free rings and their relations},
  Academic press, 1971.}

\newcommand{\CorneaRanickiRigGlu}
{O.Cornea, A.Ranicki,
\emph{
Rigidity and glueing for Morse and Novikov complexes},  e-print: AT.0107221.
Journal publication:
Journal of the European Mathematical Society 5, 343--394 (2003).}

\newcommand{\BassHellerSwanWhPoly}
{H.Bass, A.Heller, R.G.Swan,
\emph{The Whitehead group of a polynomial extension},
Inst. Hautes Etudes Sci. Publ. Math. {\bff 22} (1964),
61--79
}






\newcommand{\Conway}
{J. H. Conway,
\emph{
An enumeration of
knots and links, and some
of their algebraic properties},
Computational
Problems in Abstract Algebra,
Pergamon Press,
New York, 1970,  329--358.}

\newcommand{\ClarkHG}
{B. Clark,
\emph{
The Heegaard genus of manifolds obtained
by surgery on links and knots},
Internat. J. Math. and Math. Sci. Vol.3 No.3 (1980), 
583-589. }

\newcommand{\CorneaRanickiRG}
{O.Cornea, A.Ranicki,
\emph{
Rigidity and glueing for Morse and Novikov complexes},  e-print: AT.0107221.
Journal publ:
Journal of the European Mathematical Society 5, 343--394 (2003).}

\newcommand{\DieudonneElemAnal}
{   J.Dieudonn\'e,
\emph{El\'ements d'analyse, Vol. III
},
Gauthier-Villars, 1970.
}

\newcommand{\DieudonneFoundAnal}
{   J.Dieudonn\'e,
\emph{Foundations of modern analysis
},
Academic press, 1960.
}

\newcommand{\DollGB}
{H. Doll,
\emph{A generalized bridge number},
Math. Ann.  {\bff 294}, 1992, pp. 701--717. }

\newcommand{\DoldAT}
{   A.Dold,
\emph{Lectures on Algebraic Topology
},
Springer,  1972.
}

\newcommand{\DonaldsonTQFTSW}
{S. Donaldson,
\emph{Topological field theories and formulae of Casson and Meng-Taubes},
Geometry and Topology Monographs, Proceedings of the Kirbyfest v. 2, 1999, pp. 87 -- 102. }

\newcommand{\FarberExN}
{M.Farber
 \emph{ Exactness of Novikov inequalities },
Functionalnyi  Analiz i ego Prilozheniya  {\bff
 19},
 1985.(in Russian), Engl. Transl:
 Functional Analysis and Applications  {\bff  19},  1985.
 }
 
\newcommand{\FarberLuSch}
{M.Farber
 \emph{Lusternik - Schnirelman theory for closed 1-forms}, Commentarii Math. Helv., vol 75(2000), pages 156 -170.}
 
 \newcommand{\FarberTopCloF}
{M.Farber
 \emph{Topology of closed 1-forms and their critical points}, Topology. vol 40(2001), pages 235 - 258.}

\newcommand{\FarberTF}
{M. Farber
 \emph{Topology of closed 1-forms},
Mathematical Surveys and Monographs, 108.
American Mathematical Society, Providence, RI, 2004. }

\newcommand{\FarberRanickiMN}
{ M. Farber, A. Ranicki,
\emph{ The Morse-Novikov theory of circle-valued functions
and noncommutative localization,}
E-print:
dg-ga/9812122, Journal Publication:
Proc. 1998 Moscow Conference
for the 60th Birthday of S. P. Novikov, tr.
 Mat. Inst. Steklova, {\bf 225}, 1999,
381 -- 388.
}



\newcommand{\FarrellBAMS}
{F.T. Farrell,
\emph{The obstruction to fibering a manifold over a circle},
Bull. Amer. Math. Soc. {\bff 73} (1967),737 -- 740.
}

\newcommand{\FarrellFibInd}
{F.T. Farrell,
\emph{The obstruction to fibering a manifold over a circle},
Indiana Univ.J., {\bff 21} (1971), p. 315 -- 346.
}

\newcommand{\FarrellHsiangK}
{F.T.Farrell, W.-C.Hsiang,
\emph{A formula for $K_1R_\alpha[T]$},
Proc. Symp. Pure Math., Vol. {\bff 17} (1968), 192--218.}




\newcommand{\FelshtynNielReid}
{A.Fel'shtyn,
\emph{
Dynamical zeta functions,
Nielsen Theory and Reidemeister torsion,}
preprint: ESI 539 of The Erwin Schr\"odinger International Institute for
Mathematical Physics,
Journal Publication:
Mem. Amer. Math. Soc. 147 (2000), no. 699.
}




\newcommand{\FranksHomolDyn}
{   J.Franks
\emph{Homology and dynamical systems},
CBMS Reg. Conf. vol. 49, AMS, Providence 1982.
}




\newcommand{\FullerInd}{  F.B.Fuller,
\emph{An index of fixed point type for periodic orbits},
Amer. J.Math
{\bff 89},
(1967)
133--148.
}




\newcommand{\FriedHomolId}{  D.Fried,
\emph{Homological Identities for closed orbits},
Inv. Math. {\bff 71}, (1983) 419--442.
}




\newcommand{\FriedPerPoTwi}{  D.Fried,
\emph{Periodic points and twisted coefficients},
Lect. Notes in Math.,
{\bff 1007},
(1983)
261--293.
}




\newcommand{\FriedNewZeta}{  D.Fried,
\emph{Flow equivalence, hyperbolic systems and a new zeta function
for flows},
Comm. Math. Helv.,
{\bff 57},
(1982)
237--259.
}



\newcommand{\GallotHulinLafontaineRG}
{  S.Gallot, D.Hulin, J.Lafontaine,
\emph{Riemannian Geometry},
 Universitext.
 Springer, 2004.
}

\newcommand{\GeogheganNicasTraTors}
{  R.Geoghegan, A.Nicas,
\emph{Trace and torsion in the theory of flows},
Topology,
{\bff 33},
(1994)
683--719.
}



\newcommand{\GabaiMur}
{  D.Gabai,
\emph{The Murasugi sum is a natural geometric operation},
Contemp. Math.,
{\bff 20},
(1983)
131--143. }

\newcommand{\GeogheganNicasParamLef}
{  R.Geoghegan, A.Nicas,
\emph{Parameterized Lefschetz-Nielsen fixed
 point theory and Hochshild homology traces},
Amer. J. Math.,
{\bff 116},
(1994)
397--446. }



\newcommand{\GeogheganNicasEulChar}
{  R.Geoghegan, A.Nicas,
\emph{Higher Euler characteristics 1},
L'Enseignement Math\'ematique,
{\bff 41},
(1995)
3--62
}



\newcommand{\GrobmanHomeoDAN}
{D.M.Grobman,
\emph{Homeomorphisms of systems of differential equations},
DAN SSSR
{\bff 128}
(1959), 880-881.
\quad Engl. transl:
 Sov.Math.Dokl.
{\bff 128}
(1959), 880-881.}

\newcommand{\GrobmanTopClass}
{D.M.Grobman,
\emph{
Topological classification of the \nei~ of
a singular point in $n$-dimensional space},
Matem. Sbornik
{\bff  56}
(1962), 77-94.
}

\newcommand{\HadamardIter}
{ J.Hadamard,
\emph{Sur l'it\'eration et les solutions asymptotiques
des \'equations diff\'erentielles},
Bull. Soc. Math. France,
{\bff 29},
(1901)
224 -- 228.
}

\newcommand{\HatcherAT}
{A. Hatcher,
\emph{Algebraic topology},
Cambridge University Press, Cambridge 2002.
}




\newcommand{\GodaKitanoMorifugiReidTwi}
{ H.Goda, T.Kitano, T.Morifugi,
\emph{Reidemeister torsion,
twisted Alexander polynomial
and fibered knots }, preprint, 2002;
Journal Publication:
Comment. Math. Helv. 80(2005), 51--61.
}



\newcommand{\GodaMorifugiTwiAl}
{ H.Goda,  T.Morifugi,
\emph{
Twisted Alexander polynomial
for $SL(2,\CCC)$-representations
and fibered knots
},
preprint, 2002
}



\newcommand{\GodaHeg}
{ H.Goda,
\emph{Heegaard splitting for sutured
manifolds and Murasugi sum,
}
Osaka J.Math.
{\bff 29},
(1992)
21 -- 40.
}




\newcommand{\GodaMatsudaMorifujiFlo}
{ H.Goda, H. Matsuda, T. Morifuji,
\emph{ Knot Floer Homology of (1, 1)-Knots,
}
Geometriae Dedicata, Volume 112, Number 1, 
April 2005 , pp. 197-214.
}

\newcommand{\GodaHandNum}
{ H.Goda,
\emph{On handle number of
Seifert surfaces in $S^3$
},
Osaka J.Math.
{\bff 30},
(1993)
63 -- 80.
}



\newcommand{\GodaPajitnovTwiNov}
{H.Goda, A.Pajitnov,
\emph{ Twisted Novikov homology and 
circle-valued Morse theory for knots and links},
e-print: math.GT/0312374,
Journal Publication: Osaka Journal of Mathematics,
v.{\bff 42} No. 3,  2005, p. 557 -- 572.}

\newcommand{\GodaMatsudaPajitnovMNHeg}
{H.Goda, H. Matsuda, A.Pajitnov,
\emph{ Morse-Novikov theory, Heegaard splittings and closed orbits of gradient flows },
e-print: arXiv:0709.3153.}

\newcommand{\HamiltonInvNash}
{ R.S.Hamilton,
\emph{On the inverse
function theorem of Nash and Moser
},
Bull. AMS,
{\bff 7},
(1982)
65 -- 229.
}



\newcommand{\HartmanODE}
{ P.Hartman
\emph{Ordinary differential equations},
Wiley; New York,  1964.
}

\newcommand{\HarveyLawsonFinVol}
{R.Harvey, B.Lawson,
\emph{Morse theory and finite volume flows},
Ann.Math, 
{\bff 153}, 2001
1 - 25.
}




\newcommand{\HarveyLawsonMorStok}
{R.Harvey, B.Lawson,
\emph{Morse theory and Stokes theorem},
Surveys in Differential Geometry, 
{\bff 7},
65 -- 229.
}



 
\newcommand{\HarveyMinerviniMorNov}
{
\emph{Morse Novikov theory and cohomology with finite supports},
e-print: math.DG/0212295. 
}




\newcommand{\HigmanUnits}
{ G.Higman,
\emph{Units in group rings},
Proc. London Math. Soc.,
{\bff 46},
(1940)
231 -- 248.
}



\newcommand{\HirasawaRudolph}
{M. Hirasawa, Lee Rudolph,
\emph{Constructions of Morse maps for 
knots and links, and upper bounds on the Morse-Novikov number},
math.GT/0311134, to appear in J. Knot Theory Ramifications.
}

\newcommand{\RudolphMuS}
{Lee Rudolph,
\emph{Murasugi sums of Morse maps to the circle, 
Morse-Novikov numbers, and free genus of knots},
math.GT/0108006.
}

\newcommand{\HirschDT}
{ M.Hirsch,
\emph{Differential Topology},
Springer,  1976.
}

 \newcommand{\HuebschMorseBowl}
{W.Huebsch, M.Morse,
\emph{The bowl theorem and a model
 non-degenerate function},
Proc.Nat.Acad.Sci.U.S.A., {\bff 51}, (1964), 49 - 51.
}




\newcommand{\HusemollerFB}
{ D.Husemoller,
\emph{Fibre bundles
},
McGraw-Hill, 1966.
}

\newcommand{\HutchingsLeeCVMRTSW}
{M.Hutchings, Y.-J.Lee,
\emph{Circle-valued Morse theory, Reidemeister torsion
and Seiberg-Witten invariants of 3-manifolds},
\quad E-print:
 dg-ga/9612004,  3  Dec 1996,
 journal publication:
Topology,
{\bff 38},
(1999),
861 -- 888.
}


\newcommand{\HutchingsLeeCVMRT}
{ M.Hutchings, Y-J.Lee,
\emph{  Circle-valued Morse theory and Reidemeister torsion,}
Geometry and Topology,
{\bff 3},
(1999),
369 -- 396.
}




\newcommand{\irwinStM}
{   M.Irwin,
\emph{ On the Stable Manifold Theorem},
Bull. London Math. Soc,
{\bff 2},
(1970),
196 -- 198.
}


\newcommand{\IrwinDynSys}
{ M.C.Irwin, 
\emph{Smooth dynamical systems},
Academic Press,  1980.
}

\newcommand{\JiangPerOrb}
{  B.Jiang,
\emph{ Estimation of the
 number of periodic orbits}
Preprint:  Universit\"at Heidelberg,
Mathematisches Institut, Heft 65, Mai 1993,
Journal Publication:
Pac. J. Math. 172, No.1, 151-185 (1996).
}


\newcommand{\JiangWangTwiRep}
{  B.Jiang, S.Wang,
\emph{ Twisted topological invariants associated with representations}, in:
Topics in Knot Theory, Kluwer, 1993,
pages 211 -- 227.
}

\newcommand{\JungHoonLee}
{Jung Hoon Lee,
\emph{An upper bound for tunnel number of a knot using free genus},
Talk at 4th East Asian School of knots, 
http://faculty.ms.u-tokyo.ac.jp/$\sim$topology/EAS4slides/JungHoonLee.pdf}

\newcommand{\KanenobuPret}
{  Taizo Kanenobu,
\emph{ The augmentation subgroup of a pretzel link, }
Mathematics Seminar Notes of Kobe University, 
v. 7,  363-384 (1979).
}

\newcommand{\KawauchiDist}
{  Akio Kawauchi,
\emph{ Distance between links by zero-linking twists, }
Kobe Journal of Mathematics, 
 13, No.2, 183-190 (1996).
}

\newcommand{\KelleyGT}
{ J.L.Kelley, 
\emph{General Topology},
van NOstrand, 1957.
}

\newcommand{\KiangFP}
{   Kiang Tsai-han, 
\emph{The theory of Fixed point classes},
 Springer, 1989.
}

\newcommand{\KirkLivingstonTwiAl}
{P.Kirk, Ch. Livingston,
\emph{Twisted Alexander invariants,
 Reidemeister torsion, and
 Casson-Gordon invariants, }
Topology,
{\bff 38}
1999,
p. 635 -- 661.
}


\newcommand{\KinoshitaTerasaka}
{S. Kinoshita, H. Terasaka,
\emph{On unions of knots},
Osaka Math. J. {\bff 9} (1957), 131-153.}

\newcommand{\KitanoTwiAlTors}
{T.Kitano,
\emph{
Twisted Alexander polynomial
and  Reidemeister torsion, }
Pacific J.Mat,
{\bff 174}
(1996)
p. 431 -- 442.
}


\newcommand{\KlingenbergLectGeod}
{   W.~Klingenberg,
\emph{Lectures on closed geodesics},
 Springer, 1978.
}

\newcommand{\KobayashiRieckGRTN}
{T. Kobayashi, Y. Rieck,
\emph{
On the growth rate of tunnel number of knots}
math.GT/0402025, 
J. Reine Angew. Math. 592 (2006) 63 -- 78.
}

\newcommand{\KobayashiDTN}
{T. Kobayashi,
\emph{
A construction of arbitrarily high
degeneration of tunnel number of knots under connected sum, }
J. Knot Theory Ramifications {\bff 3}, (1994) p. 179-186.
}

\newcommand{\KodamaKNOT}
{K.Kodama,
\emph{ KNOT program},
\newline
 http://www.math.kobe-u.ac.jp/~kodama/knot.html.}

\newcommand{\KozlovskiiIntAn}
{O.S.Kozlovskii,
\emph{The dynamics
of intersections of
analytic manifolds
}, Doklady ANSSSR,
{\bf 323}
(1992), in Russian. 
\quad English translation:
 Sov.Math.Dokl.
{\bff 45}
(1992), 425--427.}

\newcommand{\KupkaContrib}
{   I.Kupka,
\emph{
Contribution \`a la th\'eorie des champs g\'en\'eriques},
Contributions to Differential equations,
{\bff 2}
    (1963    ), 457--484,
    {\bff 3}
    (1964    ), 411--420.
 }

\newcommand{\LamSerreConj}
{    T.Y.Lam,
\emph{Serre's Conjecture,         }
Lecture Notes in Mathematics {\bff 635}, (1978) 227 p.
}

\newcommand{\LangAlg}
{    S.Lang,
\emph{Algebra,}
Addison-Wesley (1965)
}

\newcommand{\Latour}
{F.Latour,
\emph{
Existence de 1-formes ferm\'ees non
singuli\`eres dans une classe de cohomologie de de Rham},
Publ. IHES {\bff 80}
(1995), 135 -- 194.
}



 
\newcommand{\LatschevThesis}
{ J.Latschev,
\emph{A generalization of the Morse complex},
PhD thesis, SUNY, 1998.}



\newcommand{\LatschevMatAnn}
{J.Latschev,
\emph{Gradient flows of Morse-Bott functions},
Mathematische Annalen,  vol. 318, no. 4, (2000), pp. 731-759.
}

\newcommand{\LaudenbachTS}
{F.Laudenbach,
\emph{On the Thom-Smale complex},~
In {\it An extension of a theorem by Cheeger and M\"uller}, ~
Asterisque, {\bff  205} (1992)
p. 219 -- 233.}

\newcommand{\LaudenbachSikoravIsoHam}
{   F.Laudenbach, J.-C.Sikorav,
\emph{
Persistance d'intersection avec la section
nulle au cours d'une isotopic hamiltonienne
 dans un fibre cotangent},
Invent.~Math.~
{\bff 82}
    (1985     ),
pp. 349--357.
 }

 \newcommand{\Lazarev}
{  A.Lazarev,
\emph{
Novikov homologies in knot theory},
Mathematical Notes 
{\bff 51}
    (1992    ),
pp. 259--261.
 }

\newcommand{\LueckUnivLef}
{ W.L\"uck,
\emph{The Universal Functorial Lefschetz Invariant }
Fundam. Math. 161, No.1-2, 167-215 (1999).
}

\newcommand{\LiapounovProbStab}
{   Liapounov A.M.,
\emph{Probl\`eme g\'en\'eral de 
la stabilit\'e du mouvement},
Thesis, Harkov, 1892 (in Russian),
\quad
French Transl:
 Ann. Fac. Sci. Toulouse, 
{\bff 9}, (1907), p. 203-474.
 \quad Engl transl:
\emph{The General Problem of the Stability of Motion},
Taylor ~$\&$~ Francis, 1992.
}

\newcommand{\MarkTTQFT}
{  T. Mark, 
\emph{ Torsion, TQFT, and Seiberg-Witten invariants of 3-manifolds},
Geom. Topol. 6 (2002), 27 - 58.
}

\newcommand{\MasseyHCT}
{   W.Massey,
\emph{ Homology and cohomology theory},
 Marcel Dekker, 1978.
}

\newcommand{\MilnorHCob}
{   J.~Milnor,
\emph{Lectures on the
h-cobordism theorem},
 Princeton University
Press,
 1965.
}

\newcommand{\MilnorCycCov}
{J.Milnor,
\emph{ Infinite cyclic coverings},
In: Conference on the topology of manifolds,
(1968). }

\newcommand{\MilnorSevenSphere}
{J.Milnor,
\emph{On manifolds homeomorphic to the 7-sphere},
Annals of Mathematics,
{\bff 64}
(1956),
399 - 405.
(1968)}

\newcommand{\MilnorKT}
{   J.~Milnor,
\emph{ Introduction to algebraic K-theory},
 Princeton University Press, 1971.
}

\newcommand{\MilnorMT}
{   J.~Milnor,
\emph{ Morse theory},
 Princeton University Press, 1963.
}

\newcommand{\MilnorWT}
{J.Milnor,
\emph{ Whitehead Torsion},
Bull. Amer. Math. Soc.
{\bff 72}
(1966),
358 - 426.
}

\newcommand{\MilnorStasheffChCl}
{J.Milnor and J.Stasheff,
\emph{ Characteristic Classes},
 Princeton University Press,
 1974.}

\newcommand{\MinerviniThesis}
{G.Minervini,
\emph{A Current Approach to  Morse and Novikov theories},
PhD thesis, Universit\`a La Sapienza, Roma, 2003.}

\newcommand{\MorimotoSATN}
{K. Morimoto,
\emph{On the super additivity of tunnel number of knots},
Math. Ann., {\bff 317}, No. 3 (2000), p. 489-508
}

\newcommand{\MorimotoATN}
{K. Morimoto,
\emph{On the  additivity of tunnel number of knots},
Topology Appl., {\bff 53}, No. 3 (1993), p. 37--66.
}

\newcommand{\MorimotoSakumaYokotaTN}
{K. Morimoto, M. Sakuma, Y. Yokota, 
\emph{Examples of tunnel number one  knots which have the property "1+1=3"},
Math. Proc. Camb. Phil. Soc., {\bff 119} (1996), p. 113--118.
}

\newcommand{\MorseCrP}
{M.Morse,
\emph{Relations between the critical points of a real function of $n$ independent variables},
Trans. Amer. Math. Soc.
{\bff 15} (1925), 345-396. }

\newcommand{\MorseCV}
{M.Morse,
\emph{Calculus of Variations in the Large},
  American Mathematical Society Colloquium Publications,
Vol.18,
 1934. }

\newcommand{\MunkresEDT}
{  J.R.Munkres,
\emph{Elementary differential topology},
Annals of Math. Studies,
Vol.54, Pinceton
 1963.}

\newcommand{\NovikovShmeltserKiEq}
{S.P.Novikov, I.Shmel'tser
\emph{Periodic solutions of the Kirchhof equations for the free motion
of a rigid body in an ideal fluid and the extended
Lyusternik-Shnirel'man-Morse 
theory, 1
}
Funkstionalnyi analiz i ego prilozheniya, 
(1981), 54--66.
\quad   Engl transl:
 Funktional analysis and applications,
{\bff 15}
(1981),
       197--207. }

\newcommand{\NovikovVarPer}
{S.P.Novikov,
  \emph{Variational methods and
  periodic soluions of equations
  of  Kirchhof  type, 2},
  Funkstionalnyi analiz i ego prilozheniya, {\bff 15}
(1981), Engl transl:
 Functional analysis and applications,
{\bff 15}
(1981), 263--274. }

 \newcommand{\NovikovMVDAN}
  {S.P.Novikov,
\emph{ Many-valued functions
 and functionals. An analogue of Morse theory  },
Doklady AN SSSR,
{\bf 260}
(1981),  31-35 (in Russian),
\quad English translation:
 Sov.Math.Dokl.
{\bff 24}
(1981), 222-226. }

 \newcommand{\NovikovRMS}
 {S.P.Novikov,
\emph{The hamiltonian formalism and a
multivalued analogue of
Morse theory,
}
Uspekhi Mat. Nauk,
{\bff 37}    (1982),  3-49(in Russian),
\quad English translation:
 Russ. Math. Surveys,
{\bff 37} (1982),
 1--56.}

\newcommand{\NovikovBH}
 {S.P. Novikov,
\emph{Bloch homology, critical points of functions and closed 1-forms}
DAN SSSR,  {\bff 287} (1982), $N^\circ$ 6,
\quad English translation:
Soviet Math. Dokl. {\bff 287} (1986), $N^\circ$ 2.
}

\newcommand{\PalisDeMelo}
{   J.Palis, Jr.,W.de Melo,
\emph{ Geometric  theory of dynamical systems},
 Springer, 1982.
}

\newcommand{\PeixotoKS}
{M.Peixoto,
\emph{ On an approximation theorem of Kupka and Smale},
J. Diff.Eq.
{\bff 3}
(1966),
423 -- 430.
}

\newcommand{\PitcherIneq}
{E.Pitcher,
\emph{ Inequalities of critical point theory},
Bull. Amer. Math. Soc.
{\bff 64}
(1958),
1-30.
}

\newcommand{\NovikovQuasiPer}
{S.P.Novikov,
\emph{Quasiperiodic Structures in topology},
in the  book:  Topological Methods in Modern Mathematics,
 Publish or Perish,
 1993,
pp. 223--235.  }

\newcommand{\PolyaSzego}
{G.Polya, G.Szeg\"o,
\emph{Aufgaben und Lehrs\"{a}tze 
aus der Analysis}
Springer, 1964}

\newcommand{\PostnikovMT}
{   M.M.Postnikov,
\emph{Introduction to Morse theory
},
 Moscow, Nauka,  1971, in Russian
}

\newcommand{\PozniakTriangMorse}
{M.Pozniak,
\emph{Triangulation of compact smooth manifolds and Morse theory}
(University of Warwick preprint, 11/1990,
published posthumously as a part of the thesis of M.Pozniak
in Translations of AMS, 2000)}

\newcommand{\PozniakMorNov}
{M.Pozniak,
\emph{The Morse complex, Novikov Cohomology and Fredholm Theory}
(University of Warwick preprint, 08/1991,
published posthumously as a part
of the thesis of M.Pozniak
in Translations of AMS, 2000)}

\newcommand{\QuillenProM}
{ D.Quillen,
\emph{Projective modules over polynomial
 rings},
Inventiones Math.,
 {\bff 36} (1976), 167-171.}

\newcommand{\RanickiFinDom}
{ A.Ranicki,
\emph{Finite domination
and  Novikov rings},
preprint: 1993, Edinburgh Yniversity, journal article:
Topology,
 {\bff 34} (1995), 619--632.}

\newcommand{\RanickiAlgNov}
{ A.Ranicki,
\emph{The algebraic construction
of the Novikov complex of a circle-valued
Morse function},

\quad  E-print: math.AT/9903090,
Journal Publ: 
Mathematische Annalen 322, 745--785 (2002). 
}

\newcommand{\RanickiMTandNH}
{ A.Ranicki,
\emph{
Circle valued Morse theory and Novikov homology },
School on High-dimensional Manifold Topology, 543--569, ICTP Trieste (2002). e-print AT.0111317.
}

\newcommand{\RanickiKL}
{   A.A.Ranicki,
\emph{Lower $K$- and $L$-theory,         }
LMS Lecture Notes 178, Cambridge, 1992
}

\newcommand{\RanickiKno}
{   A.A.Ranicki,
\emph{High-dimensional knot theory,         }
Springer, 1998
}

\newcommand{\RanickiTors}
{   A.A.Ranicki,
\emph{The algebraic theory of torsion I.,}
Lecture Notes in Mathematics {\bff 1126} (1985), 199--237.
}

\newcommand{\ReidemeisterLins}
{K.Reidemeister,
\emph{Homotopieringe und Linsenr\"aume},
Hamburger Abhandl. {\bff 11} (1938), 102--109.}

\newcommand{\Riley}
{R.Riley,
\textit{
Homomorphisms of knots on finite groups},
Mathematics of Computation, v. 25, 
{\bff 115}, 1971, p. 603 -- 619.}

\newcommand{\RockafellarCA}
{   R.T.Rockafellar,
\emph{Convex Analysis},
Princeton University Press (1970).
}

\newcommand{\RolfsenKL}
{   D.Rolfsen,
\emph{Knots and Links},
Publish or Perish (1976, 1990).
}

\newcommand{\RosenbergAlgK}
{   J.Rosenberg,
\emph{Algebraic $K$-theory
 and its applications},
Springer, (1994).
}

\newcommand{\HirasawaRudolphConsM}
{   M. Hirasawa, L.  Rudolph,
\emph{Constructions of Morse maps for knots and links, and upper bounds on the Morse-Novikov number},
e-prints mathGT/0108006 (2001), math.GT/0311134 (2003).
}

\newcommand{\ScharlemannSchultensAHTN}
{M. Scharlemenn, J. Schultens,
\emph{Annuli in generalized Heegaard splitting and degeneration of tunnel number},
Math. Ann {\bff 317} (200) No. 4, 783--820.
}

\newcommand{\Schubert}
{H. Schubert,
\emph{
\"Uber eine numerische Knoteninvariante},
Math. Z. {\bff 61} (1954), 245--288.
}

\newcommand{\SchuetzOrbs}
{D.Sch\"utz,
\emph{Gradient flows of closed 1-forms and their closed orbits},
e-print:
math.DG/0009055,
~ journal article:
Forum Math. 14(2002) 509--537.
}

\newcommand{\SchuetzOneParam}
{ D.Sch\"utz,
\emph{One-parameter fixed point theory and gradient flows
of closed 1-forms},
e-print:
math.DG/0104245,
~ journal article:
K-theory, 25(2002), 59-97.
}

\newcommand{\SharkoF}
{ V.V.Sharko,
\emph{ Funktsii na mnogoobraziyah
(algebraicheskie i topologicheskie aspekty)}, Kiev,
Naukova dumka
(1990),
 31-35.
   \quad English translation:
 }

\newcommand{\SheihamNcCohn}
{D.Sheiham
\emph{ Non-commutative Characteristic Polynomials and Cohn Localization},
e-print:   math.RA/0104158 
Journal publ: Journal of the London Mathematical Society (2) Vol 64 (2001) no.1 pp13-28.}

\newcommand{\SheihamWhLoc}
{D.Sheiham
\emph{Whitehead Groups of Localizations and the Endomorphism Class Group}, 
e-print:
  math.KT/0209311. Journal publ:
 Journal of Algebra, Vol 270 (2003) Issue 1, 261-280.}

\newcommand{\Siebenmann}
{L.Siebenmann,
\emph{A total Whitehead torsion obstruction to fibering over the circle},
Comment. Math. Helv. {\bff 45} (1970), 1--48.}

\newcommand{\SilverWilliamsAlexNTK}
{D. Silver, S. Williams,
\emph{
 Twisted Alexander polynomials detect the unknot}, Algebraic and Geometric Topology, 6 (2006), 1893-1907.
arXiv:math/0604084v3 [math.GT]}

\newcommand{\SikoravFlo}
{J.-Cl; Sikorav,
\emph{
Homologie associ\'ee
a une fonctionnelle
(d'apr\`es Floer)},
Asterisque 201 - 202 - 203, 1991,
S\'eminaire Bourbaki {\bff 733}, 
1990 -- 1991,
November 1991.}

\newcommand{\SikoravDisj}
{J.-Cl.~Sikorav,
\emph{Un probleme de disjonction par
isotopie symplectique dans un
fibr\'e cotangent},
 Ann.~Scient.~Ecole~Norm.~Sup.,{\bff 19}
 (1986),  543--552.}

\newcommand{\SikoravThesis}
{ J.-Cl.~Sikorav,
\emph{
Points fixes de diff\'eomorphismes
symplectiques, intersections de sous-vari\'et\'es
lagrangiennes, et singularit\'es de un-formes ferm\'ees
}
Th\'ese de Doctorat d'Etat Es
Sciences Math\'ematiques,
Universit\'e Paris-Sud, Centre d'Orsay, 1987}

\newcommand{\SmaleStruct}
{  S.~Smale, \emph{On the structure of manifolds},
 Am.~J.~Math., {\bff 84} (1962)  387--399.}

\newcommand{\SmaleGenPoi}
{  S.~Smale, \emph{Generalized Poincar\'e
onjecture in dimensions greater than four},
 Ann.~Math., {\bff 74} (1961)  391--406.}

\newcommand{\SmaleGradDyn}
{  S.~Smale, \emph{On gradient dynamical systems},
 Am.~J.~Math., {\bff 74} (1961)  199--206.}

\newcommand{\SmaleDyn}
{  S.~Smale,
\emph{Differential dynamical systems},
 Bull. Amer. Math. Soc. {\bff 73} (1967)
747--817.}

\newcommand{\SmaleStableM}
{  S.~Smale,
\emph{Stable \ma s for differential
equations and diffeomorphisms
},
Ann. Scuola Norm. Superiore Pisa, {\bff 18} (1963)
97--16.}

\newcommand{\SmalePoi}
{  S.~Smale,
\emph{Generalized Poincare's conjecture in dimensions
greater than four},
 Ann.~Math.,
{\bff 74} (1961)
 391--406.}

 \newcommand{\SpanierAT}
{ E.H.Spanier,
\emph{Algebraic topology},
  McGraw-Hill,
(1966)}

\newcommand{\StallingsFib}
{J.Stallings ,
\emph{On fibering certain 3-manifolds},
 Proc. 1961 Georgia conference on the Topology of 3-manifolds,
Prentice-Hall, 1962, pp. 95--100.}

\newcommand{\StongNC}
{R.Stong,
\emph{Notes on Cobordism theory},
Princeton, New Jersey, 1968}

\newcommand{\SwitzerAT}
{ R.M.Switzer,
\emph{Algebraic topology -- homotopy and homology},
  Springer,
( 1975)}

\newcommand{\ThomCell}
{ R.Thom,
\emph{Sur une partition en cellules associ\'ee
\`a une fonction sur une vari\'et\'e},
     Comptes Rendus de l'Acad\'emie de Sciences,
{\bff 228}
(1949),
 973--975.
}

\newcommand{\TuraevAlPol}
{V.Turaev,
\emph{ The Alexander polynomial 
of a three-dimensional manifold},
Math. USSR Sbornik, v. 26 (1975), No. 3, p. 313-329.
}

\newcommand{\TuraevReidTors}
{
\emph{
Reidemeister Torsion in knot theory,}
Russian Math. Surveys, 41:1 (1986), 119 - 182.
}

\newcommand{\TuraevEulStr}
{
\emph{ Euler structures, nonsingular vector fields
and torsions of Reidemeister type},
Math. USSR Izvestia 34:3 (1990), 627 - 662.
}

\newcommand{\WadaTALP}
{ M.Wada,
\emph{Twisted Alexander polynomial
for finitely presentable groups},
 Topology,
{\bff 33}
(1994),
 241 -- 256.
}

\newcommand{\VanDerWaerden}
{ B.L.van der Waerden,
\emph{Algebra 1},
Springer, 
( 1971)}

\newcommand{\WeibelHA}
{ C.A.Weibel,
\emph{An introduction to homological algebra },
Cambridge University press,
( 1997)}

\newcommand{\WellsAC}
{ R.O.Wells,
\emph{Differential Analysis on complex manifolds},
Prentice-Hall,
( 1973)}

\newcommand{\WittenSMT}
{ E.Witten,
 \emph{Supersymmetry and Morse theory},
 Journal of Diff.~Geom.,
{\bff 17} (1982)
661 -- 692.
}

\newcommand{\WittenMonoFour}
{ E.Witten,
 \emph{Monopoles and Four-Manifolds},
 Math. Res. Letters, 
{\bff 1} (1994)
769 -- 796.
}

\newcommand{\WhiteheadSHT}
{ J.H.C.Whitehead
 \emph{Simple homotopy types},
 Amer. J. Math.,
{\bff 72} (1952)
 pp. 1- 57
}


\newcommand{\PajitnovAn}
{  A.Pajitnov
 \emph{An analytic
proof of the real part of the Novikov inequalities}
 DAN SSSR,
{\bff 293},
 1987
no. 6 (in Russian);
\quad Engl transl:
 Soviet Math. Dokl.{\bff 35}  (1987), no. 2.}

\newcommand{\PajitnovShDAN}
{ A.Pajitnov,
 \emph{On the sharpness
of the inequalities of Novikov type in the case $\pi_1(M)={\bold Z}^m$
for Morse forms whose cohomology classes are in general position     },
DAN SSSR
{\bff 306},
 1989
no. 3.(In Russian)
\quad English translation:
 Soviet Math. Dokl.{\bff 39}  (1989), no. 3.}

\newcommand{\PajitnovSbor}
{  A.Pajitnov
\emph{On the sharpness of
 Novikov-type inequalities for
manifolds with free abelian fundamental group.} Matem. Sbornik,
(1989)
no. 11.
\quad English translation:
\emph{
 },
 Math. USSR Sbornik,
{\bff 68}
 (1991),
 351 - 389.
}

\newcommand{\PazhitnovModLoc}
{
A.V.Pazhitnov,
\emph{Modules over some localizations
of the ring of Laurent polynomials},
Mathematical Notes, {\bff 46}, 1989, no. 5.}

\newcommand{\PajitnovToul}
{ A.V.Pajitnov, \emph{ On the Novikov
complex for rational Morse forms},
\quad preprint:
Institut for Matematik og datalogi, Odense Universitet
Preprints 1991, No 12, Oct. 1991;
~ journal article:
Annales de la Facult\'e de Sciences de
Toulouse {\bf 4}  (1995), 297--338.
}

\newcommand{\PajitnovSurg}
{ A.V.Pajitnov,
\emph{
Surgery on the Novikov Complex},
Preprint:
Rapport de Recherche CNRS URA 758,  Nantes,   1993;
journal article:
K-theory {\bff 10} (1996),  323-412.
 }

\newcommand{\PajitnovMRL}
{  A.V.Pajitnov,
\emph{Rationality and exponential growth
properties of the boundary operators in the Novikov
Complex},
Mathematical Research Letters,
{\bff 3}
(1996),
  541-548.
 }

\newcommand{\PajitnovAsymp}
{  A.V.Pajitnov,
\emph{   On the asymptotics of
Morse numbers of finite covers of
manifold
},

\quad E-print:
math.DG/9810136, 22 Oct 1998,\quad
journal article:
 Topology,
\textbf{38}, No. 3,  pp. 529 -- 541
(1999).

}

\newcommand{\PajitnovAdv}
{  A.V.Pajitnov,
\emph{   $C^0$-generic properties of
boundary operators in the Novikov
complex },
\quad E-print:
math.DG/9812157, 29 Dec 1998,
journal article:
Advances in Mathematical Sciences,
 vol. 197, 1999, p.29 -- 117.
}

\newcommand{\PajitnovRMS}
{  A.V.Pajitnov,
\emph{ Simple homotopy type of Novikov complex
and $\zeta$-function of the gradient flow, }
\quad E-print:
dg-ga/970614 9 July 1997;
journal article:
Russian Mathematical Surveys,
\textbf{54}
(1999), 117 -- 170.}

\newcommand{\PajitnovIncRat}
{  A.V.Pajitnov,
{\it   The incidence coefficients in the Novikov
complex are
generically rational functions,}
\quad E-print:  dg-ga/9603006 14 March  96,
journal article: Algebra i Analiz,
{\bff 9}, no.5 (1997),  92--139 (in Russian),
\quad English translation:
Sankt-Petersbourg Mathematical Journal
\textbf{9}
(1998),
no. 5, p. 969 -- 1006.
}

\newcommand{\PajitnovWitt}
{  A.V.Pajitnov,
\emph{Closed orbits of gradient
flows and logarithms of non-abelian Witt vectors},
\quad E-print:
 math.DG/9908010, 2 Aug. 1999
journal article:
  K-theory, Vol. 21 No. 4, 2000, pp. 301 -- 324.
}

\newcommand{\PajitnovRanickiWN}
{A.V.Pajitnov, A.Ranicki,
\emph{The Whitehead group of the Novikov ring},
\quad E-print:
 math.AT/0012031, 5 dec 2000,
 journal article:
  K-theory, Vol. 21 No. 4, 2000.
}

\newcommand{\PajitnovGrad}
{A.V.Pajitnov,
\emph{$C^0$-topology in Morse theory},
\quad E-print:
 math.DG/0303195, 5 dec 2000.
}

\newcommand{\PajitnovWeberRudolphMNK}
{  A.V.Pajitnov, C.Weber, L.Rudolph,
\emph{ Morse-Novikov number for knots and links},~
Algebra i Analiz,
{\bff 13}, no.3 (2001),
(in Russian),
English translation:
Sankt-Petersbourg Mathematical Journal.
{\bff 13}, no.3 (2002), p. 417 -- 426.
}

\newcommand{\PajitnovCloOrb}
{  A.V.Pajitnov,
\emph{ Counting closed orbits
of gradients of circle-valued maps},~
E-print:  math.DG/0104273 28 Apr. 2001,
journal article:
Algebra i Analiz,
{\bff 14}, no.3 (2002),  92--139
(in Russian),
English translation:
Sankt-Petersbourg Mathematical Journal.
{\bff 14}, no.3 (2003).
}

\newcommand{\PajitnovCVMT}
{ A.V.Pajitnov,
\emph{Circle-valued Morse theory},
de Gruyter studies in Mathematics, vol. 32,
(2006)}

\title{On the tunnel number and the 
Morse-Novikov number of knots}
\author{A.V.Pajitnov}
\address{Laboratoire Math\'ematiques Jean Leray 
UMR 6629,
Universit\'e de Nantes,
Facult\'e des Sciences,
2, rue de la Houssini\`ere,
44072, Nantes, Cedex}
\email{ pajitnov@math.univ-nantes.fr}
\begin{abstract}
 Let $L$ be a link in $S^3$;
denote by $\MM\NN(L)$
the Morse-Novikov number of $L$
and by $t(L)$ the tunnel number of $L$.
We prove that $\MM\NN(L)\leq 2t(L)$
and deduce several corollaries.
\end{abstract}
\maketitle


\section{Introduction}
\subsection{Background}
Let $L$ be a link in $S^3$, that is, an embedding of several copies 
of $S^1$ to $S^3$. First off, we recall  the definition of three numerical
invariants of $L$. In the sequel $N(L)$ denotes a 
closed tubular \nei~ of $L$.
\pa
{\it A. ~ Tunnel Number.}  An arc $\g$ in $S^3$
is called {\it a tunnel} for $L$ if $\g\cap L$ consists of the two endpoints of $\g$. The tunnel number $t(L)$ is the minimal number $m$
of disjoint tunnels $\g_1,\ldots, \g_m$ such that 
the closure of 
$S^3\sm N(L\cup\g_1\cup\ldots \cup\g_m)$ 
is a handlebody.
The tunnel number was introduced by B. Clark in \cite{ClarkHG};
this invariant was  studied in the works of K. Morimoto,
M.Sakuma, Y. Yokota, T. Kobayashi, M. Scharlemenn, J. Schultens
and others (see \cite{MorimotoSATN},
\cite{MorimotoATN},
\cite{MorimotoSakumaYokotaTN},
\cite{ScharlemannSchultensAHTN}).

For any two knots $K_1, K_2$ we have
$t(K_1\krest K_2)\leq t(K_1)+t(K_2)+1$.
In the paper \cite{KobayashiRieckGRTN}
 T.Kobayashi and Y. Rieck 
defined the growth rate for a knot $K$ by the formula
$$
gr_t(K)=
\underset{m\to\infty}{\rm lim~ sup} \frac {t(mK)-mt(K)}{m-1}
$$
where $mK$ stands for the connected sum of $m$
 copies of the knot $K$.
We have $1\geq gr_t(K)\geq -t(K)$. 
\pa
{\it B. Bridge numbers.}
Let $S^3=H_1\cup H_2$ be a Heegaard decomposition of $S^3$;
put $\Sigma = H_1\cap H_2$, and $g=g(\Sigma)$.
We say (following H. Doll \cite{DollGB}) that $L$ is in a $n$-bridge position
\wrt~ $\Sigma$ if $\Sigma$ intersects $L$ in $2n$ points
and $\Sigma\cap H_i$ is a union of $n$ trivial arcs in $H_i$
for $i=1,2$.
The $g$-bridge number $b_g(L)$ of $L$ is defiined as the minimal number $n$ such that $L$ can be put in a $n$-bridge position \wrt~ a Heegaard decomposition of genus $g$
(thus $b_0(L)$ is the classical bridge number as defined n the paper \cite{Schubert}  of H. Schubert).
We have
$$
t(L)\leq g+b_g(L)-1.
$$
\pa
{\it C. Morse-Novikov numbers.}
A {\it framing} of $L$ is a diffeomorphism
$\phi:L\times D^2\to N(L)$.
Let $C_L$ denote the closure of $S^3\sm N(L)$.
A Morse function $f:C_L\to S^1$ is called {\it regular}
if its restriction to the boundary $\pr N(L)$ satisfies 
the following relation:
$(f\circ \phi)(l,z)=\frac z{|z|}$.
A regular Morse 
function has finite number of critical points;
 the number of the critical points of $f$ of index $i$ wil be denoted by $m_i(f)$; the total number of critical points of $f$
will be denoted by $m(f)$.
The minimal value of $m(f)$ over all possible framings $\phi$ and all
possible Morse maps $f:C_L\to S^1$
is called {\it the Morse-Novikov number of the link $L$}
and denoted by $\MM\NN(L)$
(see \cite{PajitnovWeberRudolphMNK}).
The Morse-Novikov theory implies that
$$
\MM\NN(L)\geq 2(b_1(L)+q_1(L))
$$
where $b_1(L)$ and $q_1(L)$ are the  {\it Novikov numbers}
defined as follows.
Let $\bar C_L$ be the infinite cyclic covering
induced by $f$ from the covering $\RRR\to S^1$.
Denote the ring $\ZZZ[t,t^{-1}]$ by $\L$, and
the ring $\ZZZ((t))$ by $\wh{\L}$.
Then $b_1(L)$ and $q_1(L)$ are respectively the 
rank and torsion numbers of the module
$H_1(\bar C_L)\tens{\L}\wh \L$.
In case when the Novikov numbers are not sufficient
 to determine the $\MM\NN(L)$ the twisted Novikov numbers
(introduced by H. Goda and the author in 
\cite{GodaPajitnovTwiNov})
are useful.

As for the upper 
bounds for $\MM\NN(L)$ 
not much is known.
M. Hirasawa proved that for every 2-bridge knot $K$
we have $\MM\NN(K)\leq 2$ (unpublished). In the papers
\cite{RudolphMuS} and 
\cite{HirasawaRudolph}
of Lee Rudolph and M. Hirasawa  it is proved that
$\MM\NN(K)\leq 4g_f(K)$
where $g_f(K)$ is the {\it free genus } of $K$, that is,
the minimal possible genus of a Seifert surface $\Sigma$
bounding $K$ such that $S^3\sm \Sigma$ is an open handlebody.

\subsection{Main results}
\label{s:main_res}

The  main result of this work is 
\beth\label{t:mn_t}
For every link $L$ in $S^3$ we have
\begin{equation}\label{f:mn_t}
\MM\NN(L)\leq 2t(L).
\end{equation}
\end{theo}

The following corollaries are easily deduced.
\beco
For every $g$ we have
$$
\MM\NN(L)\leq 2(g+b_g(L)-1).
$$
\end{coro}
\beco
For every (1,1)-knot $K$ we have 
$\MM\NN(K)\leq 2$.
\end{coro}
\beco
For every link $L$ we have
$$
q_1(L)+b_1(L)\leq t(L).
$$
\end{coro}
\beco
For every knot $K$ 
$$gr_t(K)\geq -t(K)+q_1(K).$$
\end{coro}

\section{Proof of Theorem \ref{t:mn_t}}

Let $m=t(L)$.
Pick a framing $\phi:L\times D^2\to N(L)$.
 Then the manifold $C_L=\overline{S^3\sm N(L)}$ is obtained from $\pr C_L$ by attaching $m$ one-handles and then attaching 
a handlebody of genus $(m+1)$ 
to the resulting cobordism. 
Thus we obtain a Morse function $g:C_L\to\RRR$ which is constant on $\pr C_L$ and has the following Morse numbers:
$m_0(g)=0,\ m_1(g)=m,\ m_2(g)=m+1,\ m_3(g)=1$.
Pick any Morse map $h:C_L\to S^1$ such that $h|\pr C_L$ is the canonical fibration:
$(h\circ\phi)(l,z)=\frac z{|z|}$.
Consider a closed 1-form $\o_\e=dg+\e dh$.
For $\e>0$ sufficiently small $\o_\e$ is a Morse form with the same Morse numbers as $dg$. Therefore 
the form 
$$\frac 1\e\o_\e=d\Big(\frac 1\e g+h\Big)$$
is the differential of a Morse map $g_1:C_L\to S^1$
having the required behaviour on $\pr C_L$.
The map $g_1$ has one local maximum, and the standard elimination procedure (see for example \cite{PajitnovWeberRudolphMNK} for details)
gives us a Morse function
$f:C_L\to S^1$ with
$m_0(f)=0,\ m_1(f)=m,\ m_2(f)=m,\ m_3(f)=0$.
Thus $\MM\NN(L)\leq 2m$.

\section{Examples, and further remarks}
\label{s:remarks}

A theorem of M. Hirasawa says that 
$\MM\NN(K)\leq 2$ if $K$ is a two-bridge knot.
Since $t(K)\leq b(K)-1$ our theorem implies 
this result. Observe that the proof of the M. Hirasawa's theorem 
uses the H. Schubert's classification of 
2-bridge knots, and can not be generalized to the 
case of arbitrary bridge number.

The inequality \rrf{f:mn_t}
implies also the upper bound 
$$
\MM\NN(K)\leq 4g_f(K)
$$
obtained by Lee Rudolph and M. Hirasawa
(see \cite{RudolphMuS}, 
\cite{HirasawaRudolph}).
Indeed  Jung Hoon Lee 
\cite{JungHoonLee}
has shown that 
$t(K)\leq 2g_f(K)$.

In many cases the estimate of Theorem \ref{t:mn_t}
is better than the free genus estimate.
For example, for a pretzel knot 
$K=P(-2, m,n)$ where $m,n\geq 3 $ are odd numbers,
we have $g(K)=\frac {m+n}2$
(see \cite{GodaMatsudaMorifujiFlo}), so that 
$g_f(K)\geq \frac {m+n}2$. On the other hand
$t(K)=1$.

\section{The tunnel number and the 
homology with local coefficients
}
\label{s:tunn_morse}

Let $L$ be a link in $S^3$, put $m=t(L)$.
As we have observed  in the previous section
there is a Morse function $g:C_L\to\RRR$
such that $g$ is constant on $\pr C_L$ and takes there
its minimal value, and with the Morse numbers as follows:
$m_0(g)=0,\ m_1(g)=m,\ m_2(g)=m+1,\ m_3(g)=1$.
The function $-g$ provides a handle decomposition of
the manifold $C_L$, with $(m+1)$ one-handles,
therefore we have the usual homological estimate $m+1\geq \mu_\ZZZ(H_1(C_L,\ZZZ))$.
\footnote{\ \ 
For a finitely generated $R$-module $T$ we denote 
by  $\mu_R(T)$  the minimal number of generators of $T$.}
For the case of knots this estimate is trivial,
however in some cases we can improve it  using  
homology with local coefficients.
Let $\rho:\pi_1(C_L)\to GL(q,R)$ be a 
{\it right } representation
(that is, $\rho(ab)=\rho(b)\rho(a)$ for every $a,b$).
Denote by $\wi C_L$ the universal covering of $C_L$.
The homology of the chain complex 
$$
C_*(\wi C_L,\rho)=R^q\tens{\rho}C_*(\wi C_L)
$$
is called {\it the homology with local coefficients $\rho$}
or {\it $\r$-twisted homology}
and denoted by 
$H_*(C_L,\rho)$.
If $R$ is the principal ideal domain, then we have
\begin{equation}\lb{f:tun_homloc}
m+1\geq \frac 1q\ \mu_R\Big(H_1(C_L,\rho)\Big).
\end{equation}
In what follows we will concentrate on the case of
knots. For a knot $K$ consider a meridional embedding
$i:S^1\to C_K$. Given a right representation 
$\r:\pi_1(C_K)\to GL(q,R)$ we can induce 
it to $S^1$ by $i$ and obtain a local coefficient
system $i^*\rho$ on $S^1$.

The following proposition is an easy corollary of the 
main theorem of the paper of D. Silver and S. Williams
\cite{SilverWilliamsAlexNTK}.

\bepr\label{p:nontriv_locsys}
Let $K$ be any knot in $S^3$. 
Then there is a right representation 
$\g:\pi_1(C_K)\to GL(q,R)$ with $R$ a principal ideal domain
such that 
\been\item[(i)] $H_1(C_K,\g)\not=0$,
\item[(ii)] $H_1(S^1, i^*\g)=0$.
\enen
\end{prop}
\Prf
Let us first recall briefly the Silver-Williams theorem.
Consider the meridional 
homomorphism $\xi:\pi_1(C_K)\to\ZZZ$
as a homomorphism of $\pi_1(C_K)$ to $\L^\bu=GL(1,\L)$,
where $\L=\ZZZ[t,t^{-1}]$.
For a right representation $\t:\pi_1(C_K)\to GL(q,\ZZZ)$
form the tensor product $\r=\xi\otimes\t:\pi_1(C_K)\to GL(q,\L)$.
Consider the $\L$-module $\gB=H_1(C_K,\r)$
and choose a free resolution for $\gB$: 
\begin{equation}\label{f:reso}
\xymatrix{0 & \ar[l] \gB & \ar[l] {\L^r} & \ar[l]_p {\L^k}
& \ar[l] {\ldots}}
\end{equation}
where $k\geq r$. The GCD of the ideal of $\L$
generated by the $r\times r$-minors of $p$ is called the 
{\it twisted Alexander polynomial } of $k$ \wrt~ $\r$;
we will denote it 
$\D(K,\t)$ (it is defined up to multiplication by $\pm t^i$).
The Silver-Williams theorem says that for every $K$ there is 
a representation $\t$ such that $\D(K,\t)$ is not a unit 
of $\L$.
Pick such a representation $\t$ and consider two cases:
\pa
{ 1) }  
$\D(K,\t)$ is a monomial, that is  $\D(K,\t)=at^n$
with $a\in\ZZZ, \ a\not=\pm 1$. In this case 
define the representation $\g=\wh\r$
to be the composition of $\r$ with the natural inclusion
$GL(q,\L)\sbs GL(q,\wh\L)$.
The  $\wh\r$-twisted homology is 
$$H_1(C_K,\wh\rho)=H_1(C_K,\rho)\tens{\L}\wh\L.$$
We can obtain a free $\wh\L$-resolution for 
this module by tensoring
the resolution \rrf{f:reso} by $\wh\L$ over $\L$.
Since the GCD of elements of $\L$ remains the same when we 
extend the ring $\L$ to $\wh \L$ (see \cite{PajitnovWeberRudolphMNK},
Lemma 2.3), the GCD of the $r\times r$-minors of the matrix $\wh p$
equals $a$. Since $\wh\L$ is a principal ideal domain 
we deduce that $H_1(C_K,\wh\rho)$ is non-zero 
and moreover it contains a 
cyclic direct summand.
The property (ii) is easy to check.

\pa
2) $\D(K,\t)$ is a polynomial of non-zero degree.
In this case consider the ring $\L_\QQQ=\QQQ[t,t^{-1}]$.
This ring is principal and $\D(K,\t)$ is not invertible in
it; define the representation $\g=\wi\r$
to be the composition of $\r$ with the natural inclusion
$GL(q,\L)\sbs GL(q,\L_\QQQ)$.
The same argument as for the point 1) works here 
as well.
$\qs$

\bepr\label{p:tungr_lin}
Let $K$ be any knot in $S^3$. Then there is $\l>0$
such that for every $n\in\NNN$ we have
$t(nK)\geq n\l -1$.
\enpr
\Prf
Pick  a representation $\g:\pi_1(C_K)\to GL(q,R)$
satisfying the conclusion of 
Proposition \ref{p:nontriv_locsys}. 
The module $H_1(C_{K},\r)$ contains then
a cyclic $R$-submodule $T$.
\bele\label{l:lem_tun}
For any $n\geq 1$ there is a right representation
$\g_n:\pi_1(C_{nK})\to GL(q,R)$ such that the module
$\gB_n=H_1(C_{nK},\g_n)$ contains a submodule isomorphic 
to $nT$.
\end{lemm}
\Prf We proceed by induction in $n$.
Denote by $\mu\in \pi_1(C_K)$
the meridional element.
Assume that we have constructed $\g_n:\pi_1(C_{nK})\to GL(q,R)$
in such a way that $\g_n(\mu)=\g(\mu)$.
The group $\pi_1(C_{(n+1)K})$ is isomorphic to the 
amalgamated  product of the groups 
$\pi_1(C_K)$ and $\pi_1(C_{nK})$ over the subgroup $\ZZZ$
included to both groups via the embedding of the meridian.
Let $\g_{n+1}:\pi_1(C_{(n+1)K})\to GL(q,R)$
be the product of the representations
$\g$ and $\g_n$.  Using the property (ii) from 
\ref{p:nontriv_locsys} and the Mayer-Vietoris exact sequence
it is easy to deduce that the module
$\gB_{(n+1)}=H_1(C_{(n+1)K},\g_{(n+1)})$
contains a submodule isomorphic to the direct sum
of $T$ and $nT$. $\qs$

The previous Lemma implies that
$
\mu_R(\gB_n)\geq n. 
$
Our proposition follows, 
since
$$
t(nK)+1\geq \frac 1q \
\mu_R(\gB_n)\geq \frac nq.
$$
$\qs$
\beco\label{c:gr_lin}
For any knot $K$ we have
$$
gr_t(K)> -t(K).
$$
$\qs$
\end{coro}

\section{Generalizations and a question}
\label{s:gen_open}

The generalization
of the results of the Section \ref{s:main_res}
to the case of knots and links in an arbitrary closed 3-manifold
are straightforward. The same goes for the formula
\rrf{f:tun_homloc}. On the other hand  it is not clear at 
all whether the Proposition \ref{p:nontriv_locsys} and 
\ref{p:tungr_lin} admit such generalizations, since 
the analogs of Silver-Williams theorem for arbitrary 
three-manifolds seem to be out of reach for the moment.

{\bf Question.} Are the inequalities \rrf{f:tun_homloc}
sufficient to determine the tunnel number for every link?
In other words is it true that
\begin{equation}
t(L)+1= \max_{\r}\Big( \frac 1q\ \mu_R\big(H_1(C_L,\rho)\big)\Big).
\end{equation}
where $\r$ ranges over all right representations
$\pi_1(C_L)\to GL(q,R)$?
\section{Acknowledgements}
\label{s:ackn}

This work was completed during my stay 
as a GCOE visitor at 
the Graduate School of Mathematical Sciences
at Tokyo University.
Many thanks to T. Kohno for warm hospitality.

I am grateful to Hiroshi Goda for 
many valuable discussions on knot theory 
since several years.

\end{document}